\documentclass[preprint,12pt,3p]{elsarticle}

\usepackage{amsmath, amssymb, amsfonts, amsbsy, amsthm, latexsym, epsfig, txfonts, palatino, pifont, tikz, color, graphicx}
\usetikzlibrary {positioning}
\definecolor {processblue}{cmyk}{0.96,0,0,0}

\theoremstyle{definition}

\newtheorem{df}{Definition}

\journal{Mathematical Biosciences}

\makeatletter
\g@addto@macro\@floatboxreset\centering
\makeatother

\begin{document}

\begin{frontmatter}

\title{Local Immunodeficiency: Minimal Networks and Stability}

\author{Leonid Bunimovich}
\ead{bunimovh@math.gatech.edu}

\author{Longmei Shu\corref{cor1}}
\address{School of Mathematics,
Georgia Institute of Technology,
Atlanta, GA 30332-0160 USA}

\cortext[cor1]{Corresponding author}

\ead{lshu6@math.gatech.edu}

\begin{abstract}
Some basic aspects of the recently discovered phenomenon of local immunodeficiency \cite{pnas} generated by antigenic cooperation in cross-immunoreactivity (CR) networks are investigated. We prove that local immunodeficiency (LI) that's stable under perturbations already occurs in very small networks and under general conditions on their parameters. Therefore our results are applicable not only to Hepatitis C where CR networks are known to be large \cite{pnas}, but also to other diseases with CR. A major necessary feature of such networks is the non-homogeneity of their topology. It is also shown that one can construct larger CR networks with stable LI by using small networks with stable LI as their building blocks. Our results imply that stable LI occurs in networks with quite general topology. In particular, the scale-free property of a CR network, assumed in \cite{pnas}, is not required.
\end{abstract}


\end{frontmatter}



\section{Introduction}

Cross-immunoreactivity (CR) is a well known phenomenon which was observed in the studies of AIDS, influenza, Hepatitis C, dengue and other diseases (see e.g. \cite{hattori,campo,nowak00,yoshioka,wodarz,nowak90,nowak91,nowak91s}). In a  nutshell, CR means that the generation of antibodies to some antigen (virus) can be stimulated by other antigens. Therefore CR generates (indirectly, i.e. via the corresponding antibodies) interactions between the antigens. For a long time CR was recognized as an important phenomenon in the in-host dynamics of various diseases and was used in building their mathematical models \cite{nowak90,nowak91,wodarz,nowak00}.

However, in all these models CR was incorporated as a mean-field process where all interactions between different antigens (viruses) are assumed to have the same strength. Recent experiments with Hepatitis C viruses demonstrated that this assumption is incorrect, and instead the CR network (CRN) has a very complicated structure (topology) which resembles the topology of scale-free networks \cite{campo,hattori}.

A new model for the dynamics of Hepatitis C (HC) \cite{pnas} is conceptually simpler than the previous ones (see e.g. \cite{wodarz}). In fact the new model involves only two (necessary) types of variables, the population sizes of various types of viruses and the population sizes of their corresponding antibodies, in immunological models. For instance, the HC model in \cite{wodarz} contains three more types of variables, namely the population sizes of infected and of non-infected hepatocytes as well as a total (mean field) CR response.

This fact naturally causes some doubts and suspicion. Indeed, how can a simpler model have richer dynamics? The reason is that our model is just conceptually simpler,  it actually contains more parameters. Different pairs of viruses generally have different strengths of interaction in the CRN, but in the old model they were all equal to each other.

Traditionally, to describe new experimental findings which old models fail to reproduce, one makes a more complicated mathematical model by adding more variables or more equations. The model introduced in \cite{pnas} is based on new specially conducted experiments \cite{campo,hattori} which proved essential heterogeneity of the CRN.
Although the model in \cite{pnas} was dealing with dynamics of HC, it provides a model of evolution for any disease which has cross-immunoreactivity.
The paper \cite{pnas} analyzed the dynamics of this new model numerically. Scale-free CRNs of sizes 500-1000 were generated and numerical simulations were performed on them.

The main result was the discovery of a new phenomenon \cite{pnas}, Local Immunodeficiency (LI), which showed up in all of the several hundred simulations.
Namely, in all these simulations, the pool of HC viruses got partitioned into three types. The first type consists of \emph{persistent} viruses that have large population sizes but virtually zero immune response against them. In other words, \emph{persistent} viruses remain undetected by the human immune system. Thus a clear immunodeficiency (with respect to persistent viruses) is present. It is called \cite{pnas} \emph{local immunodeficiency} because it is completely determined by the localized positions of the \emph{persistent} viruses in the CRN. Observe, however, that generally it may happen that only specific types of antigens are "qualified" to be persistent viruses. Only special biological experiments may clarify this issue. 

\emph{Persistent} viruses enjoy such a relaxing life because the second type, \emph{altruistic} viruses, sacrifice themselves to protect the \emph{persistent} viruses from the immune system. Concentrations of \emph{altruistic} viruses are very small but they carry almost the entire immune response against all of the in-host population of viruses. Again, we need further experimental biological studies to determine which antigens can and which can not play a role of altruistic viruses.  The rest (third type) of viruses plays a much smaller role in the HC evolution \cite{pnas}. In what follows we call these viruses \emph{neutral}.

In the present paper we demonstrate rigorously that local immunodeficiency is a much more general phenomenon than one may conclude from the results of \cite{pnas}.

First, we prove that stable LI already appears in a specific CRN with only three nodes under general conditions. These conditions are expressed as realistic inequalities between parameters of the model. Therefore LI is likely to appear in all diseases with cross-immunoreactivity. Indeed, because of a very high mutation rate of HC viruses in host, the corresponding CRNs are very large \cite{pnas}.  Since both small and large CRNs can generate LI, one is tempted to believe that this phenomenon should be universal for all diseases with cross-immunoreactivity.

It is proved that LI is a stable state of evolution of the model in only one (out of many possible topologies) of the networks with three types of viruses, while in all two-node CRNs LI is unstable. This three-node network with stable LI is characterized by the maximal asymmetry of its structure among all networks of size three. Here by "maximal asymmetry" we mean that all the nodes have different indegrees. In this network there is one persistent node and one altruistic node while the third node is neutral.

We also prove that there are no two-node CRNs with stable LI. It should be mentioned that the two-node network with stable LI found in \cite{pnas} assumes very restrictive relations between parameters of the model, which have the form of exact equalities. Clearly such strict constraints cannot be maintained in real life situations. Indeed, only inequalities remain true under small changes of parameters, which always occur because of fluctuations of real environments. In the present paper we demonstrate that the regions in the parameter space where stable LI exists have the same dimension as the the parameter space of the model. However, it happens only in certain networks with at least three elements (types of viruses). Once again, these networks must also be sufficiently non-homogeneous, which is (qualitatively) consistent with numerical results in \cite{pnas} for large CRNs.

We then demonstrate how one can build larger CRNs with stable LI by attaching the three-node minimal network with stable LI. For instance, we proved that by combining two such networks one gets a network with five nodes where two types of viruses are persistent and two are altruistic. And the dynamics of HC with such a CRN is stable and robust. Our results were mostly obtained by direct computations. For large networks one would need numerical simulations although our rigorous results about smaller CRNs basically give a proof of concept that stable and robust LI is present in all larger networks with sufficiently non-homogeneous topology.

To justify it even more we also prove the presence of stable and robust LI in a network with seemingly mild non-homogeneity of its topology. It is important to mention that among CRNs with four nodes there are quite a few with more non-homogeneous topology than the one we studied. Therefore our results essentially prove that stable and robust LI must also be present in those CRNs. It is for this purpose that we studied a less non-homogeneous network. The proof of stable and robust LI (essentially by long direct computations) in this CRN is given in the Appendix.

It is important to mention that in this paper we are dealing with \emph{strong} LI, which is a stronger property than the one found in \cite{pnas}. Namely, we say that a certain type of virus causes strong local immunodeficiency if the immune response against it is identically zero, so completely absent. Analogously, we say that some kind of virus is altruistic if it is not present at all (i.e. its concentration is zero) but immune response against this non-existing virus is present (strictly positive).

In \cite{pnas}, instead of these identical zeroes, some (sufficiently) small quantities were considered. We call this case a \emph{weak} LI.  Clearly a \emph{weak} LI is a more general phenomenon than \emph{strong} LI. Indeed, if the strong LI takes place then the weak LI is automatically present. Thus our results imply that \emph{weak} LI does exist and is stable, under \emph{even weaker} conditions than our conditions on the existence and stability of \emph{strong} LI. Therefore it is present in an even larger variety of CRNs.

These rigorously proven results demonstrate that stable LI does not require a special scale-free structure of the CRN. In fact, it is enough that the CRN is sufficiently non-homogeneous. It is natural to expect that this condition is satisfied in real life situations because there is no reason for CRNs to be homogeneous. Non-homogeneity of CRNs is a mild and very general condition, and thus the phenomenon of local immunodeficiency should be ubiquitous for diseases with cross-immunoreactivity. 

We also show that LI is a robust phenomenon. Recall that a state of a system is stable if small variations of initial conditions result in small variations of this state, i.e. a new (perturbed) orbit stays close to the initial (unperturbed) state. On the other hand, a state of a system is robust if small variations of the system parameters (i.e. transitions to formally different systems) result in a stable state which is close to the state of the initial (unperturbed) system.  

Our results demonstrate once again that altruistic viruses, which have very small concentration but occupy central positions in the CRNs with the largest indegrees \cite{pnas}, play a key role in LI. Namely the altruistic viruses were present in all CRNs where we found stable and robust LI. All CRNs with fixed points with LI but without altruistic viruses turned out to be non-robust, i.e. the LI could be destroyed by arbitrarily small variations of parameters. This means that such cases are non-typical, i.e. they have a positive codimension (or zero volume) in the space of all systems we study. Therefore they cannot be seen in real life situations. 
(In other words, it is a zero probability event to encounter an LI without altruistic viruses.) This observation also explains why altruistic viruses were always present in the several hundred numerical experiments conducted in \cite{pnas}.
Therefore these altruistic local hubs of CRNs must be the primary targets of prevention and elimination of the corresponding diseases. This is yet another question for the future studies, both biological and computational. From a general biomedical point of view a main challenge is to understand which types of viruses could play a role of altruistic and which persistent ones. 

The structure of the paper is as follows. Section \ref{model} introduces the model. Section \ref{analysis} is devoted to a general analysis of the stability of dynamics of this model. Section \ref{size2} analyzes two-node networks. Three-node networks are studied in section \ref{size3}. Section \ref{altruistic} proves the necessity of altruistic viruses. The building of larger networks with stable LI is considered in section \ref{coex}. Lastly section \ref{dis} contains some concluding remarks. Some long technical computations are placed in the Appendix. We also put some long computations with a four-node CRN in the Appendix to demonstrate that LI appears in networks with a relatively mild non-homogeneity of their topology.

\section{Model of evolution of a disease with heterogeneous CRN}\label{model}
In this section we define the model of the HC evolution introduced in \cite{pnas}. It is important to stress again that this model is applicable to any disease with cross-immunoreactivity.

Consider any immunological model, a population of $n$ viral antigenic variants $x_i$ inducing $n$ immune responses $r_i$ in the form of antibodies (Abs). The viral variants exhibit CR which results in a CR network.  The latter is a directed weighted graph $G_{CRN}=(V,E)$, with vertices corresponding to viral variants and directed edges connecting CR variants.
Because not all interactions with Ab lead to neutralization, we consider two sets of weight functions for the CRN. These functions are defined by immune neutralization and immune stimulation matrices $U=(u_{ij})_{i,j=1}^n$ and $V=(v_{ij})_{i,j=1}^n$, where $0\le u_{ij},v_{ij}\le1$; $u_{ij}$ represents the binding affinity of Ab to $j~(r_j)$ with the $i$-th variant; and $v_{ij}$ reflects the strength of stimulation of Ab to $j~(r_j)$ by the $i$-th variant. The immune response $r_i$ against variant $x_i$ is neutralizing; i.e., $u_{ii}=v_{ii}=1$. The evolution of the antigen (virus) and antibody populations is given by the following system of ordinary differential equations (ODEs):

\begin{equation}\label{population}
\begin{split}
\dot x_i=f_ix_i-px_i\sum_{j=1}^nu_{ji}r_j,\quad i=1,\dots,n,\\
\dot r_i=c\sum_{j=1}^nx_j\frac{v_{ji}r_i}{\sum_{k=1}^nv_{jk}r_k}-br_i,\quad i=1,\dots,n.
\end{split}
\end{equation}

The viral variant $x_i$ replicates at the rate $f_i$ and is eliminated by the immune responses $r_j$ at the rates $pu_{ji}r_j$. The immune responses $r_i$ are stimulated by the $j$-th variant at the rates $cg_{ji}x_j$, where $g_{ji}=\frac{v_{ji}r_i}{\sum_{k=1}^nv_{jk}r_k}$ represents the probability of stimulation of the immune response $r_i$ by the variant $x_j$. This model (as in \cite{pnas}) allows us to incorporate the phenomenon of the original antigenic sin  \cite{tom,pan,shin,hvi,jin,alotofauthors}, which states that $x_i$ preferentially stimulates preexisting immune responses capable of binding to $x_i$. The immune response $r_i$ decays at rate $b$ in the absence of stimulation.

Here we consider the situation where the immune stimulation and neutralization coefficients are equal to constants $\alpha$ and $\beta$, respectively. To be more specific, both the immune neutralization and stimulation matrices are completely defined by the structure of the CRN, i.e.,
$$U=\text{Id}+\beta A^T,V=\text{Id}+\alpha A,$$
where $A$ is the adjacency matrix of $G_{CRN}$. In the absence of CR among viral variants the system  reduces to the model developed in \cite{nowak00} for heterogeneous viral population.
Because the neutralization of an antigen may require more than one antibodies, we assume that $0<\beta=\alpha^k<\alpha<1$ \cite{pnas}. It is important to mention that we analyze a more general model here than the one studied in \cite{pnas}, where it was assumed that all viruses replicate with the same rate. 

\section{Stationary states of the model}\label{analysis}

Fixed points of system \eqref{population} are determined by the relations 
\begin{equation}\label{stationary}
\begin{split}
f_ix_i=px_i\sum_{j=1}^nu_{ji}r_j,\quad i=1,\dots,n,\\
cr_i\sum_{j=1}^n\frac{v_{ji}x_j}{\sum_{k=1}^nv_{jk}r_k}=br_i,\quad i=1,\dots,n.
\end{split}
\end{equation}

Clearly we are interested only in such fixed points where all variables assume non-negative values, and the populations of all viruses and antibodies can not be simultaneously equal to zero.

Consider the following sets
\begin{align*}
N=\{i\in\mathbb{N},1\le i\le n\}, I=\{i\in N: x_i>0\}, J=\{i\in N: r_i>0\}.
\end{align*}

\begin{df}
We say that \emph{strong} local immunodeficiency occurs when there exists $i$ such that $x_i>0,r_i=0$, or when $P:=I\setminus J\neq\emptyset$.
\end{df}
In what follows we will call neutral nodes with $x_i=r_i=0$ the neutral idle nodes since they don't contribute to the dynamics of the network. We also will call neutral nodes with $x_i>0,r_i>0$ the neutral active nodes.
In the paper \cite{pnas} a weaker LI condition was considered. Namely a new
phenomenon of \emph{antigenic cooperation} was discovered when some (altruistic) viral variants sacrifice themselves, being strongly exposed to an immune response, for the benefit of other (persistent) viral variants which become practically hidden from the immune system. In \cite{pnas} LI was considered to be present when persistent viruses increase their population but the immune response against them was relatively small. These conditions are more practical for computer simulations, since it could take a very long time to completely eliminate some virus, but they are not very precise. Here we consider a stronger but well defined case, \emph{strong} LI. Since a strong LI automatically implies weak LI, showing that strong LI is ubiquitous for non-homogeneous CRNs demonstrates that weak LI is even more common for such networks.    

By making use of the notations introduced above we get a simpler formula for the fixed points:
\begin{align*}
\sum_{j\in N}u_{ji}r_j=r_i+\beta\sum_{ij\in E}r_j=f_i/p,\forall i\in I,\\
\sum_{j\in N}\frac{v_{ji}x_j}{\sum_{k\in J}v_{jk}r_k}=\delta_ix_i+\alpha\sum_{ji\in E}\delta_jx_j=b/c,\forall i\in J,\\
\delta_i=\frac{1}{r_i+\alpha\sum_{ik\in E}r_k}.
\end{align*}

In our parameter space $\{f_1,f_2,\dots,f_n>0,p,c,b>0,1>\alpha,\beta>0\}$, any relation having a form of equality (e.g. $f_1=\beta f_2$) defines a subset of co-dimension 1, (i.e. a non-typical subset), in the phase space of all systems described by ODE \eqref{population}. Therefore with respect to a natural phase volume such subsets have volume (measure) zero. It is practically impossible that these very restrictive conditions will be met in a real system evolving according to model \eqref{population}. Because of this we are only interested in stationary points which exist without extra conditions or under conditions expressed as inequalities between the parameters of the model.
This should be contrasted with \cite{pnas} where LI was shown to exist under much more restrictive conditions with some exact equalities between the system's parameters.

Suppose that the matrices $V=(\text{Id}+\alpha A)$ and $U=(\text{Id}+\beta A^T)$ are invertible. Denote $F=(f_1,\dots,f_n)^T$. Then one stationary point is defined by the following relation
\begin{align*}
R^*=\frac{1}{p}(U^T)^{-1}F,X^*=\frac{b}{c}(V^T)^{-1}(VR^*)=:Xr(R^*).
\end{align*}
Notice that $U,V$ are constant matrices determined by the CRN and parameters, and $F$ is a constant vector of parameters. Because of that, $R^*$ here is a constant vector, which represents the population of the antibodies. For this $R^*$, we also have a corresponding constant vector for the population of the viruses $X^*$, given as a function of $R^*$, which is denoted as $Xr$ here for convenience.

More generally, we have a stationary space defined by the following relations
$$R=R^*+\ker(U^T_I),X=Xr(R)+\ker(V^T_J),$$
where $$\ker(U^T_I)=\{w\in\mathbb{R}^n:(U^Tw)_i=0,\forall i\in I\},\ker(V^T_J)=\{w\in\mathbb{R}^n:(V^Tw)_j=0,\forall j\in J\}.$$

To verify the stability of a stationary point, we need to consider the Jacobian matrix of the right hand side of \eqref{population}. It can be written in block form as
$$J=\begin{pmatrix} A_J & B \\ C& D \end{pmatrix},$$
where
\begin{align*}
    A_J=\text{diag}(f_i-p\sum_{j=1}^nu_{ji}r_j),B_{i,j}=-px_iu_{ji},\\
    C_{i,j}=c\frac{v_{ji}r_i}{\sum_{k=1}^nv_{jk}r_k},D_{i,l}=-cr_i\sum_{j=1}^n\frac{v_{ji}x_jv_{jl}}{(\sum_{k=1}^nv_{jk}r_k)^2},l\neq i,\\
    D_{i,i}=c\sum_{j=1}^n\frac{v_{ji}x_j}{\sum_{k=1}^nv_{jk}r_k}-b-cr_i\sum_{j=1}^n\frac{v_{ji}^2x_j}{(\sum_{k=1}^nv_{jk}r_k)^2}.
\end{align*}

\section{Analysis of size 2 CRN}\label{size2}

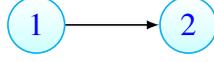
\begin {figure}[h]
\begin {tikzpicture}[-latex ,auto ,node distance =2 cm and 2cm ,on grid ,
semithick ,
state/.style ={ circle ,top color =white , bottom color = processblue!20 ,
draw,processblue , text=blue , minimum width =0.5 cm}]
\node[state] (A) {$1$};
\node[state] (B) [right =of A] {$2$};
\path (A) edge [] node[above] {} (B);
\end{tikzpicture}
\caption{size 2 CRN}
\label{2fig}
\end{figure}

We analyze the asymmetric network of size 2 (Fig.~\ref{2fig}) in this section. We consider the only asymmetric network in hope of finding LI, based on the understanding that LI requires some level of non-homogeneity of the network.

The  equations describing the evolution of these two types of viruses and antibodies are
$$\begin{cases}
\dot x_1=f_1x_1-px_1(r_1+\beta r_2),\\
\dot x_2=f_2x_2-px_2r_2,\\
\dot r_1=cx_1\frac{r_1}{r_1+\alpha r_2}-br_1,\\
\dot r_2=c(x_1\frac{\alpha r_2}{r_1+\alpha r_2}+x_2)-br_2.
\end{cases}$$

Here there is only one fixed point of interest, the one where the values of the variables are non-negative and the strong LI is present without exact equality conditions on the parameters. This fixed point is given by the relations
$$x_1=\frac{bf_1}{cp\beta},x_2=0,r_1=0,r_2=\frac{f_1}{p\beta}.$$

The Jacobian of the system is 
$$J=\begin{pmatrix}
f_1-p(r_1+\beta r_2) & 0 & -px_1 & -p\beta x_1\\
0 & f_2-pr_2 & 0 & -px_2\\
\frac{cr_1}{r_1+\alpha r_2} & 0 & \frac{cx_1\alpha r_2}{(r_1+\alpha r_2)^2}-b & -\frac{cx_1\alpha r_1}{(r_1+\alpha r_2)^2}\\
\frac{c\alpha r_2}{r_1+\alpha r_2} & c & -\frac{cx_1\alpha r_2}{(r_1+\alpha r_2)^2} & \frac{cx_1\alpha r_1}{(r_1+\alpha r_2)^2}-b
\end{pmatrix}.$$

At the fixed point the Jacobian equals

\begin{align*}
J=\begin{pmatrix}
0 & 0 & -\frac{bf_1}{c\beta} & -\frac{b}{c}f_1\\
0 & f_2-\frac{f_1}{\beta} & 0 & 0\\
0 & 0 & \frac{b}{\alpha}-b & 0\\
c & c & -\frac{b}{\alpha} & -b
\end{pmatrix}.
\end{align*}

It has the eigenvalue $\lambda=\frac{b}{\alpha}-b>0$ , and therefore this fixed point is unstable. 

It is important to mention that a stable LI for this two node network was found in \cite{pnas}. However, as we already mentioned before it has been done under unrealistic conditions. 
One can also check that the symmetric network of size 2 doesn't have a stable  LI. Detailed computations for it is listed in \ref{2sym}. Our analysis proves that no two-node network can have a stable and robust state of LI.

\section{Analysis of size 3 CRNs}\label{size3}

In this section we study the stability of dynamics of CRNs with three elements. In some of these networks there is no stable LI because of their symmetry or not enough non-homogeneity. Actually only one topology of a CRN with three elements demonstrates a stable strong LI. We present here the analysis of this size 3 CRN as well as of another one. Some other CRN is analyzed in \ref{3nodeap}.

\begin {figure}[h]
\begin {tikzpicture}[-latex ,auto ,node distance =02cm and 2cm ,on grid ,
semithick ,
state/.style ={ circle ,top color =white , bottom color = processblue!20 ,
draw,processblue , text=blue , minimum width =0.5 cm}]
\node[state] (A)
{$1$};
\node[state] (B) [right=of A] {$2$};
\node[state] (C) [right =of B] {$3$};
\path (A) edge [] node[below =0.15 cm] {} (B);
\path (B) edge [] node[below =0.15 cm] {} (C);
\end{tikzpicture}
\caption{chain-branch CRN}
\label{cb}
\end{figure}
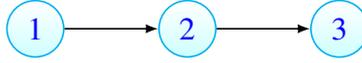

Consider at first the chain-branch CRN (Fig.~\ref{cb}). Such a network was briefly mentioned in \cite{pnas} to demonstrate that long distance action in networks may lead to LI. No studies of stability were conducted in that paper though. Also recall that here we are after \emph{robust} conditions of stable LI which would not be violated under variations of parameters. The latter always occurs because of permanently changing environments. Besides, any mathematical model (including \eqref{population} of course) is just an approximation to reality. Therefore robustness is a necessary condition for any predictive model of a real system or phenomenon.

Here system \eqref{population} becomes 
$$\begin{cases}
\dot x_1=f_1x_(1-px_1(r_1+\beta r_2),\\
\dot x_2=f_2x_2-px_2(r_2+\beta r_3),\\
\dot x_3=f_3x_3-px_3r_3,\\
\dot r_1=cx_1\frac{r_1}{r_1+\alpha r_2}-br_1,\\
\dot r_2=c(x_1\frac{\alpha r_2}{r_1+\alpha r_2}+x_2\frac{r_2}{r_2+\alpha r_3})-br_2,\\
\dot r_3=c(x_2\frac{\alpha r_3}{r_2+\alpha r_3}+x_3)-br_3.
\end{cases}$$

The fixed points with local immunodeficiency are:
$$x_1=\frac{bf_1}{cp\beta},x_2=0,x_3=0,r_1=0,r_2=\frac{f_1}{p\beta},r_3=0;$$
$$x_1=\frac{bf_1}{cp\beta},x_2=0,x_3=\frac{bf_3}{cp},r_1=0,r_2=\frac{f_1}{p\beta},r_3=\frac{f_3}{p};$$
$$x_1=\frac{bf_1}{cp},x_2=\frac{bf_2}{cp\beta},x_3=0,r_1=\frac{f_1}{p},r_2=0,r_3=\frac{f_2}{p\beta}.$$

The Jacobian of the system becomes 
$$J=\begin{pmatrix}
f_1-p(r_1+\beta r_2) & 0 & 0 & -px_1 & -p\beta x_1 & 0\\
0 & f_2-p(r_2+\beta r_3) & 0 & 0 & -px_2 & -p\beta x_2\\
0 & 0 & f_3-pr_3 & 0 & 0 & -px_3\\
\frac{cr_1}{r_1+\alpha r_2} & 0 & 0 & \frac{c\alpha x_1r_2}{(r_1+\alpha r_2)^2}-b & -\frac{c\alpha x_1r_1}{(r_1+\alpha r_2)^2} & 0\\
\frac{c\alpha r_2}{r_1+\alpha r_2} & \frac{cr_2}{r_2+\alpha r_3} & 0 & -\frac{c\alpha x_1r_2}{(r_1+\alpha r_2)^2} & \frac{c\alpha x_1r_1}{(r_1+\alpha r_2)^2}+\frac{c\alpha x_2r_3}{(r_2+\alpha r_3)^2}-b & -\frac{c\alpha x_2r_2}{(r_2+\alpha r_3)^2}\\
0 & \frac{c\alpha r_3}{r_2+\alpha r_3} & c & 0 & -\frac{c\alpha x_2r_3}{(r_2+\alpha r_3)^2} & \frac{c\alpha x_2r_2}{(r_2+\alpha r_3)^2}-b
\end{pmatrix}.$$
At the fixed point $x_1=\frac{bf_1}{cp\beta},x_2=x_3=0,r_2=\frac{f_1}{p\beta},r_1=r_3=0$, the Jacobian is 
\begin{align*}
J=\begin{pmatrix}
0 & 0 & 0 & -\frac{bf_1}{c\beta} & -\frac{b}{c}f_1 & 0\\
0 & f_2-\frac{f_1}{\beta} & 0 & 0 & 0 & 0\\
0 & 0 & f_3 & 0 & 0 & 0\\
0 & 0 & 0 & \frac{b}{\alpha}-b & 0 & 0\\
c & c & 0 & -\frac{b}{\alpha} & -b & 0\\
0 & 0 & c & 0 & 0 & -b
\end{pmatrix}.
\end{align*}
There are eigenvalues 
$\lambda=f_3,\frac{b}{\alpha}-b>0$.
Therefore this fixed point is unstable.

At the second fixed point $x_1=\frac{bf_1}{cp\beta},x_3=\frac{bf_3}{cp},x_2=0,r_2=\frac{f_1}{p\beta},r_3=\frac{f_3}{p},r_1=0$, we have
\begin{align*}
J=\begin{pmatrix}
0 & 0 & 0 & -\frac{bf_1}{c\beta} & -\frac{b}{c}f_1 & 0\\
0 & f_2-\frac{f_1}{\beta}-\beta f_3 & 0 & 0 & 0 & 0\\
0 & 0 & 0 & 0 & 0 & -\frac{b}{c}f_3\\
0 & 0 & 0 & \frac{b}{\alpha}-b & 0 & 0\\
c & \frac{cf_1}{f_1+\alpha\beta f_3} & 0 & -\frac{b}{\alpha} & -b & 0\\
0 & \frac{c\alpha\beta f_3}{f_1+\alpha\beta f_3} & c & 0 & 0 & -b
\end{pmatrix},
\end{align*}

Here $\lambda=\frac{b}{\alpha}-b>0$ is an eigenvalue, and this fixed point is also unstable.

At the fixed point $x_1=\frac{bf_1}{cp},x_2=\frac{bf_2}{cp\beta},x_3=0,r_1=\frac{f_1}{p},r_3=\frac{f_2}{p\beta},r_2=0$, the Jacobian takes the form
\begin{align*}
J=\begin{pmatrix}
0 & 0 & 0 & -\frac{b}{c}f_1 & -\frac{b}{c}\beta f_1 & 0\\
0 & 0 & 0 & 0 & -\frac{b}{c\beta}f_2 & -\frac{b}{c}f_2\\
0 & 0 & f_3-\frac{f_2}{\beta} & 0 & 0 & 0\\
c & 0 & 0 & -b & -\alpha b & 0\\
0 & 0 & 0 & 0 & \frac{b}{\alpha}+\alpha b-b & 0\\
0 & c & c & 0 & -\frac{b}{\alpha} & -b
\end{pmatrix},
\end{align*}
One eigenvalue equals
$\lambda=\frac{b}{\alpha}+\alpha b-b>0$, and hence this critical point is unstable as well. 

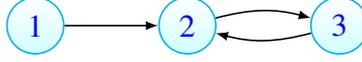
\begin {figure}[h]
\begin {tikzpicture}[-latex ,auto ,node distance =2 cm and 2cm ,on grid ,
semithick ,
state/.style ={ circle ,top color =white , bottom color = processblue!20 ,
draw,processblue , text=blue , minimum width =0.5 cm}]
\node[state] (A)
{$1$};
\node[state] (B) [right=of A] {$2$};
\node[state] (C) [right =of B] {$3$};
\path (A) edge [] node[below =0.15 cm] {} (B);
\path (B) edge [bend left =15] node[above] {} (C);
\path (C) edge [bend left =15] node[below] {} (B);
\end{tikzpicture}
\caption{branch-cycle CRN}
\label{minimal}
\end{figure}

Next we consider a CRN with three elements which has maximal asymmetry among all thirteen topologically different networks of three elements. Indeed only in this network indegrees of all three nodes are different and equal 0,2 and 1 respectively. 
In view of its essential asymmetry this network would most likely maintain LI out of all thirteen. It happened to be the case.
This network is depicted in Fig.~\ref{minimal} and we call it a branch-cycle network.

Clearly one gets a network with similar properties by relabeling the vertex 3 as 1 and vice versa.
The equations for population evolution in this case are
$$\begin{cases}
\dot x_1=f_1x_1-px_1(r_1+\beta r_2),\\
\dot x_2=f_2x_2-px_2(r_2+\beta r_3),\\
\dot x_3=f_3x_3-px_3(\beta r_2+r_3),\\
\dot r_1=cx_1\frac{r_1}{r_1+\alpha r_2}-br_1,\\
\dot r_2=c(x_1\frac{\alpha r_2}{r_1+\alpha r_2}+x_2\frac{r_2}{r_2+\alpha r_3}+x_3\frac{\alpha r_2}{\alpha r_2+r_3})-br_2,\\
\dot r_3=c(x_2\frac{\alpha r_3}{r_2+\alpha r_3}+x_3\frac{r_3}{\alpha r_2+r_3})-br_3.
\end{cases}$$

The fixed points of interest (i.e. all population sizes are non-negative, there is a strong LI, and the relations between system parameters are inequalities rather than equalities) in this case are
$$x_1=0,x_2=0,x_3=\frac{bf_3}{cp\beta},r_1=0,r_2=\frac{f_3}{p\beta},r_3=0;$$
$$x_1=\frac{bf_1}{cp\beta},x_2=0,x_3=0,r_1=0,r_2=\frac{f_1}{p\beta},r_3=0;$$
$$f_3>f_1,x_1=\frac{bf_1}{cp\beta}(1-\alpha),x_2=0,x_3=\frac{b}{cp}(f_3-f_1+\frac{\alpha}{\beta}f_1),r_1=0,r_2=\frac{f_1}{p\beta},r_3=\frac{f_3-f_1}{p};$$
$$f_3<f_1,x_1=\frac{b}{cp}(f_1-f_3+\frac{\alpha}{\beta}f_3),x_2=0,x_3=\frac{bf_3}{cp\beta}(1-\alpha),r_1=\frac{f_1-f_3}{p},r_2=\frac{f_3}{p\beta},r_3=0;$$
$$x_1=\frac{bf_1}{cp},x_2=\frac{bf_2}{cp\beta},x_3=0,r_1=\frac{f_1}{p},r_2=0,r_3=\frac{f_2}{p\beta}.$$

The Jacobian of the system is
\begin{align*}
J=\begin{pmatrix}f_1-p(r_1+\beta r_2) & 0 & 0 & -px_1 & -p\beta x_1 & 0\\
0 & f_2-p(r_2+\beta r_3) & 0 & 0 & -px_2 & -p\beta x_2\\
0 & 0 & f_3-p(\beta r_2+r_3) & 0 & -p\beta x_3 & -px_3\\
\frac{cr_1}{r_1+\alpha r_2} & 0 & 0 & \frac{cx_1\alpha r_2}{(r_1+\alpha r_2)^2}-b & -\frac{cx_1\alpha r_1}{(r_1+\alpha r_2)^2} & 0\\
\frac{c\alpha r_2}{r_1+\alpha r_2} & \frac{cr_2}{r_2+\alpha r_3} & \frac{c\alpha r_2}{\alpha r_2+r_3} & -\frac{cx_1\alpha r_2}{(r_1+\alpha r_2)^2} & A-b & -B\\
0 & \frac{c\alpha r_3}{r_2+\alpha r_3} & \frac{cr_3}{\alpha r_2+r_3} & 0 & -\frac{cx_2\alpha r_3}{(r_2+\alpha r_3)^2}-\frac{cx_3\alpha r_3}{(\alpha r_2+r_3)^2} & B-b\end{pmatrix},\\
\text{where }A=\frac{cx_1\alpha r_1}{(r_1+\alpha r_2)^2}+\frac{cx_2\alpha r_3}{(r_2+\alpha r_3)^2}+\frac{cx_3\alpha r_3}{(\alpha r_2+r_3)^2},B=\frac{cx_2\alpha r_2}{(r_2+\alpha r_3)^2}+\frac{cx_3\alpha r_2}{(\alpha r_2+r_3)^2}.
\end{align*}

At the fixed point $x_3=\frac{bf_3}{cp\beta},x_1=x_2=0,r_2=\frac{f_3}{p\beta},r_1=r_3=0$, we have 
\begin{align*}
A=0,B=b/\alpha,J=\begin{pmatrix}f_1-f_3 & 0 & 0 & 0 & 0 & 0\\
0 & f_2-\frac{f_3}{\beta} & 0 & 0 & 0 & 0\\
0 & 0 & 0 & 0 & -\frac{b}{c}f_3 & -\frac{b}{c\beta}f_3\\
0 & 0 & 0 & -b & 0 & 0\\
c & c & c & 0 & -b & -\frac{b}{\alpha}\\
0 & 0 & 0 & 0 & 0 & \frac{b}{\alpha}-b\end{pmatrix}.
\end{align*}
Because $\lambda=\frac{b}{\alpha}-b>0$ is an eigenvalue, this fixed point is unstable.

At the next fixed point $x_1=\frac{bf_1}{cp\beta},x_2=x_3=0,r_2=\frac{f_1}{p\beta},r_1=r_3=0$, we get
\begin{align*}
A=B=0,J=\begin{pmatrix}0 & 0 & 0 & -\frac{bf_1}{c\beta} & -\frac{b}{c}f_1 & 0\\
0 & f_2-\frac{f_1}{\beta} & 0 & 0 & 0 & 0\\
0 & 0 & f_3-f_1 & 0 & 0 & 0\\
0 & 0 & 0 & \frac{b}{\alpha}-b & 0 & 0\\
c & c & c & -\frac{b}{\alpha} & -b & 0\\
0 & 0 & 0 & 0 & 0 & -b
\end{pmatrix}.
\end{align*}
Hence $\lambda=\frac{b}{\alpha}-b>0$ is an eigenvalue, and this fixed point is unstable.

At the fixed point $x_1=\frac{bf_1}{cp},x_2=\frac{bf_2}{cp\beta},x_3=0,r_1=\frac{f_1}{p},r_3=\frac{f_2}{p\beta},r_2=0$, we obtain 
\begin{align*}
A=\alpha b+\frac{b}{\alpha},B=0,J=\begin{pmatrix}0 & 0 & 0 & -\frac{b}{c}f_1 & -\frac{b}{c}\beta f_1 & 0\\
0 & 0 & 0 & 0 & -\frac{b}{c\beta}f_2 & -\frac{b}{c}f_2\\
0 & 0 & f_3-\frac{f_2}{\beta} & 0 & 0 & 0\\
c & 0 & 0 & -b & -\alpha b & 0\\
0 & 0 & 0 & 0 & \alpha b+\frac{b}{\alpha}-b & 0\\
0 & c & c & 0 & -\frac{b}{\alpha} & -b
\end{pmatrix}.
\end{align*}
Then $\lambda=\alpha b+\frac{b}{\alpha}-b>0$ is an eigenvalue. This fixed point is also unstable.

For the fixed point $f_3>f_1,x_1=\frac{bf_1}{cp\beta}(1-\alpha),x_3=\frac{b}{cp}(f_3-f_1+\frac{\alpha}{\beta}f_1),x_2=0,r_2=\frac{f_1}{p\beta},r_3=\frac{f_3-f_1}{p},r_1=0$, we have 
\begin{align*}
A=\alpha b\frac{f_3-f_1}{f_3-f_1+\alpha/\beta f_1},B=b\frac{\alpha/\beta f_1}{f_3-f_1+\alpha/\beta f_1},\\
J=\begin{pmatrix}0 & 0 & 0 & -\frac{b}{c\beta}f_1(1-\alpha) & -\frac{b}{c}f_1(1-\alpha) & 0\\
0 & f_2-\frac{f_1}{\beta}-\beta(f_3-f_1) & 0 & 0 & 0 & 0\\
0 & 0 & 0 & 0 & -\frac{b\beta}{c}(f_3-f_1+\frac{\alpha}{\beta}f_1) & -\frac{b}{c}(f_3-f_1+\frac{\alpha}{\beta}f_1)\\
0 & 0 & 0 & \frac{b}{\alpha}-2b & 0 & 0\\
c & c\frac{f_1}{f_1+\alpha\beta(f_3-f_1)} & c\frac{\alpha/\beta f_1}{f_3-f_1+\alpha/\beta f_1} & b-\frac{b}{\alpha} & A-b & -B\\
0 & c\frac{\alpha\beta(f_3-f_1)}{f_1+\alpha\beta(f_3-f_1)} & c\frac{f_3-f_1}{f_3-f_1+\alpha/\beta f_1} & 0 & -A & B-b
\end{pmatrix}.
\end{align*}
Let $D=f_3-f_1+\alpha/\beta f_1,\lambda_1=f_2-f_1/\beta-\beta(f_3-f_1),\lambda_2=b/\alpha-2b$. Then
\begin{align*}
\det(\lambda I-J)=(\lambda-\lambda_1)(\lambda-\lambda_2)P(\lambda),\\
P(\lambda)=bf_1(1-\alpha)[\lambda^2+(b-B)\lambda+\frac{AD}{\alpha}]+\\\lambda\{b\beta D(\lambda+b-B)-AbD+(\lambda+b)[\lambda^2+(b-B-A)\lambda+\frac{AD}{\alpha}(1-\beta)]\}\\
=\lambda^4+b(1+\frac{(1-\alpha)(f_3-f_1)}{f_3-f_1+\alpha/\beta f_1})\lambda^3+(bf_3+b^2\frac{(1-\alpha)(f_3-f_1)}{f_3-f_1+\alpha/\beta f_1})\lambda^2\\
+b^2(1-\alpha)(f_3-f_1)(1+\frac{f_1}{f_3-f_1+\alpha/\beta f_1})\lambda+b^2(1-\alpha)f_1(f_3-f_1).
\end{align*}
One can check that all coefficients of $P(\lambda)$ are positive. It implies that $P(\lambda)$ does not have real positive roots. So in this case a stable LI is possible.
We list below a few exact values of the system parameters where stable LI is present. In each such numerical example we pick the values of the parameters to satisfy the conditions (inequalities) of existence and stability of the corresponding fixed point, and close to the literature ranges (e.g. \cite{nowak00}, \cite{pnas} and references therein).
This hand pick approach seems to be reasonable for demonstration as well as for applications. In fact in biomedical studies some parameters could be measured while the others are picked from some reasonable (accepted) ranges.
\begin{enumerate}
\item $f_1=1,f_2=3,f_3=4,b=1,\alpha=2/3,\beta=4/9$, we have $\lambda_1=-7/12<0,\lambda_2=-1/2<0$, $P(\lambda)$ has 2 pairs of conjugate complex roots, both with negative real part.
\item $f_1=1/4,f_2=1/2,f_3=1/2,b=2,\alpha=3/4,\beta=9/16$, we have $\lambda_1=-49/576<0,\lambda_2=-4/3<0$, $P(\lambda)$ has 1 pair of conjugate complex roots with negative real part and 2 distinct negative real roots.
\end{enumerate}

It is easy to see that the roots of $P(\lambda)$ depend continuously on the parameters. Therefore the set
of parameters for which the roots are real negative, or complex with negative real parts have strictly positive volume in the parameter space of the system. Thus LI in this system remains a stable type of behavior under variations of the system's parameters.

At the fixed point $f_3<f_1,x_1=\frac{b}{cp}(f_1-f_3+\frac{\alpha}{\beta}f_3),x_3=\frac{bf_3}{cp\beta}(1-\alpha),x_2=0,r_1=\frac{f_1-f_3}{p},r_2=\frac{f_3}{p\beta},r_3=0$, we have
\begin{align*}
A=\alpha b\frac{f_1-f_3}{f_1-f_3+\alpha/\beta f_3},B=\frac{b}{\alpha}-b,\\
J=\begin{pmatrix}0 & 0 & 0 & -\frac{b}{c}(f_1-f_3+\frac{\alpha}{\beta}f_3) & -\frac{b\beta}{c}(f_1-f_3+\frac{\alpha}{\beta}f_3) & 0\\
0 & f_2-\frac{f_3}{\beta} & 0 & 0 & 0 & 0\\
0 & 0 & 0 & 0 & -\frac{b}{c}f_3(1-\alpha) & -\frac{b}{c\beta}f_3(1-\alpha)\\
\frac{c(f_1-f_3)}{f_1-f_3+\alpha/\beta f_3} & 0 & 0 & \frac{b\alpha/\beta f_3}{f_1-f_3+\alpha/\beta f_3}-b & -\frac{b\alpha(f_1-f_3)}{f_1-f_3+\alpha/\beta f_3} & 0\\
\frac{c\alpha/\beta f_3}{f_1-f_3+\alpha/\beta f_3} & c & c & -\frac{b\alpha/\beta f_3}{f_1-f_3+\alpha/\beta f_3} & A-b & -B\\
0 & 0 & 0 & 0 & 0 & B-b
\end{pmatrix}.
\end{align*}
Let $D=f_1-f_3+\alpha/\beta f_3,\lambda_1=f_2-f_3/\beta,\lambda_2=b/\alpha-2b$.
Then
\begin{align*}
\det(\lambda I-J)=(\lambda-\lambda_1)(\lambda-\lambda_2)P(\lambda),\\
P(\lambda)=bf_3(1-\alpha)(\lambda^2+\frac{A}{\alpha}\lambda+\frac{AD}{\alpha})+\\
\lambda\{b\beta D(\lambda+\frac{A}{\alpha})-bAD+(\lambda+b)[\lambda^2+(\frac{A}{\alpha}-A)\lambda+\frac{AD}{\alpha}(1-\beta)]\}\\
=\lambda^4+b(1+\frac{(1-\alpha)(f_1-f_3)}{f_1-f_3+\alpha/\beta f_3})\lambda^3+(bf_1+b^2\frac{(1-\alpha)(f_1-f_3)}{f_1-f_3+\alpha/\beta f_3})\lambda^2\\
+b^2(1-\alpha)(f_1-f_3)(1+\frac{f_3}{f_1-f_3+\alpha/\beta f_3})\lambda+b^2(1-\alpha)f_3(f_1-f_3).
\end{align*}
At this point we also have that all coefficients of the polynomial $P(\lambda)$ are positive.

Again we list below several numerical values for parameters of the model where stable local immunodeficiency occurs.
\begin{enumerate}
\item $f_1=4,f_2=2,f_3=1,b=1,\alpha=2/3,\beta=4/9$, we get $\lambda_1=-1/4<0,\lambda_2=-1/2<0$. $P(\lambda)$ here has 2 pairs of complex conjugate roots, both with negative real part.
\item $f_1=1/2,f_2=1/4,f_3=1/4,b=2,\alpha=3/4,\beta=9/16$, then $\lambda_1=-7/36<0,\lambda_2=-4/3<0$. $P(\lambda)$ has 1 pair of complex conjugate roots with negative real part, and 2 distinct negative real roots.
\end{enumerate}
It follows by continuity that  there are positive volume sets in the parameter space of the model where there is a stable (i.e. practically observable) fixed point with strong local immunodeficiency.

The last size 3 CRN we consider is a 3-cycle with no stable LI. The corresponding computations are given in \ref{3nodeap}.

\section{Necessity of altruistic nodes}\label{altruistic}
We will now address a problem, whether altruistic nodes must be present in all cases of LI.

We considered all the fixed points for CRNs of sizes two and three (see \ref{listapp}). They can be separated into four groups.
\begin{itemize}
    \item A: fixed points with LI and with no extra condition on the parameters.
    \item B: fixed points with LI with conditions on the parameters in the form of inequalities.
    \item C: fixed points with LI with conditions on the parameters that involve at least one equality.
    \item D: fixed points with no LI.
\end{itemize}
One can check that fixed points in groups A and B all have altruistic nodes, while fixed points with no altruistic nodes all belong to groups C and D. So altruistic viruses are not necessary for the existence of fixed points with LI in the group C. 
However conditions on parameters in the form of equalities single out a subset of zero volume in the space of all systems we consider (when parameters in \eqref{population} assume any reasonable/permissible values). By reasonable/permissible we mean such values of parameters that make sense. For instance, negative growth rates are not permissible. 

Next we consider the existence of altruistic viruses in CRNs of arbitrary (finite) size $n$. We exclude neutral idle nodes with $x_i=r_i=0$ since they don't contribute to the dynamics. For any fixed point, assume that there are no altruistic nodes. Then  $x_i>0,\forall i=1,\dots, n$. 
This results in the following relation
\begin{equation}
    U^TR=F/p\label{alt}
\end{equation}
where $R=(r_1,\dots,r_n)^T,F=(f_1,\dots,f_n)^T$. 
It is easy to see that for \eqref{alt} to have a solution, $F$ must be in the column space of $U^T=(I+\beta A^T)^T=I+\beta A$. 

Consider now two cases.
\renewcommand{\theenumi}{\roman{enumi}}
\begin{enumerate}
\item If $U^T$ is invertible, then the column space of $U^T$ is $\mathbb{R}^n$. $F$ is always in the column space of $U^T$;
\item If $U^T$ is not invertible, then its column space is a subspace of $\mathbb{R}^n$ with a positive codimension. In other words, the condition on the parameters $f_i$'s in this case is a zero volume subset of the parameter space.
\end{enumerate}
For a fixed point to have LI, we need the vector $R$ to have at least one zero component. These vectors are on the axes and axes planes in $\mathbb{R}^n$, or the complement of the set where every component is nonzero. Hence, this is a zero volume set. 
Consider again two cases.
\begin{enumerate}
    \item If $U^T$ is invertible, then $F=pU^TR$ is also on a zero volume set.
    \item If $U^T$ is not invertible, then \eqref{alt} has either none or infinitely many solutions. Therefore if $R$ has a solution, it has one solution where some component is zero. However in the previous step we already showed that if $U^T$ is not invertible, $F$ must belong to a zero measure subspace.
\end{enumerate}

In conclusion, formally altruistic viruses are not necessary for the existence of LI. But the conditions on the parameters for fixed points to have persistent nodes without altruistic nodes are only satisfied on a zero measure subset of the parameter space. Therefore, practically speaking, altruistic viruses form a necessary component of local immunodeficiency.

\section{Building larger networks with stable \& robust LI}\label{coex}
\begin {figure}[h]
\begin {tikzpicture}[-latex ,auto ,node distance =2 cm and 2cm ,on grid ,
semithick ,
state/.style ={ circle ,top color =white , bottom color = processblue!20 ,
draw,processblue , text=blue , minimum width =0.5 cm}]
\node[state] (A)
{$1$};
\node[state] (B) [right=of A] {$2$};
\node[state] (C) [right =of B] {$3$};
\node[state] (D) [left =of A] {$4$};
\node[state] (E) [left =of D] {$5$};
\path (A) edge [] node[below =0.15 cm] {} (B);
\path (B) edge [bend left =15] node[above] {} (C);
\path (C) edge [bend left =15] node[below] {} (B);
\path (A) edge [] node[below=0.15cm] {}(D);
\path (E) edge[bend left =15] node[above] {} (D);
\path (D) edge [bend left=15] node[below] {} (E);
\end{tikzpicture}
\caption{size 5 CRN}
\label{5node}
\end{figure}
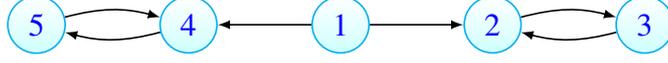

In this section we demonstrate how one can construct CRNs with multiple nodes with LI. In other words, we construct a CRN with several persistent nodes which remain hidden from the host's immune system because they are protected by the altruistic viruses. To do this we put together two identical size 3 CRNs with stable LI found in section \ref{size3}. We prove that the corresponding size 5 CRN has a fixed point with two persistent nodes and two altruistic nodes. We also demonstrate the stability of strong LI for this specific state. 
Consider the following network in Fig.~\ref{5node}.

The model \eqref{population} equations for this network are

$$\begin{cases}
\dot x_1=f_1x_1-px_1(r_1+\beta r_2+\beta r_4),\\
\dot x_2=f_2x_2-px_2(r_2+\beta r_3),\\
\dot x_3=f_3x_3-px_3(\beta r_2+r_3),\\
\dot x_4=f_4x_4-px_4(r_4+\beta r_5),\\
\dot x_5=f_5x_5-px_5(\beta r_4+r_5),\\
\dot r_1=cx_1\frac{r_1}{r_1+\alpha r_2+\alpha r_4}-br_1,\\
\dot r_2=c(x_1\frac{\alpha r_2}{r_1+\alpha r_2+\alpha r_4}+x_2\frac{r_2}{r_2+\alpha r_3}+x_3\frac{\alpha r_2}{\alpha r_2+r_3})-br_2,\\
\dot r_3=c(x_2\frac{\alpha r_3}{r_2+\alpha r_3}+x_3\frac{r_3}{\alpha r_2+r_3})-br_3,\\
\dot r_4=c(x_1\frac{\alpha r_4}{r_1+\alpha r_2+\alpha r_4}+x_4\frac{r_4}{r_4+\alpha r_5}+x_5\frac{\alpha r_4}{\alpha r_4+r_5})-br_4,\\
\dot r_5=c(x_4\frac{\alpha r_5}{r_4+\alpha r_5}+x_5\frac{r_5}{\alpha r_4+r_5})-br_5.
\end{cases}$$

Here we mirrored the chain-branch network about node 1. We are not going to try to compute all possible fixed points with LI this time. In general, based on a vague rule (there is always an arrow going from the persistent node to the altruistic node, and the altruistic node typically has a high indegree), one can make a guess and pick a node to be altruistic and another to be persistent. Then a specific fixed node with LI can be computed based on the guess through a relatively straightforward process. However, finding all possible fixed points with LI is more complicated. In the current 5-node CRN, we want LIs at both ends of this network, in the form of $x_5>0,r_5=0,x_4=0,r_4>0,x_1>0,r_1>0,x_2=0,r_2>0,x_3>0,r_3=0$. The corresponding fixed point is 
$$f_1-f_3-f_5>0,x_1=\frac{b}{cp}(f_1-f_3-f_5+\frac{\alpha}{\beta}f_3+\frac{\alpha}{\beta}f_5),r_1=\frac{f_1-f_3-f_5}{p},$$
$$x_2=0,r_2=\frac{f_3}{p\beta},x_3=\frac{bf_3}{cp\beta}(1-\alpha),r_3=0,x_4=0,r_4=\frac{f_5}{p\beta},x_5=\frac{bf_5}{cp\beta}(1-\alpha),r_5=0.$$

The Jacobian is $$J=\begin{pmatrix}A & B\\C & D\end{pmatrix},$$
$$A=\begin{pmatrix}0 & 0 & 0 & 0 & 0\\
0 & f_2-pr_2 & 0 & 0 & 0\\
0 & 0 & 0 & 0 & 0\\
0 & 0 & 0 & f_4-pr_4 & 0\\
0 & 0 & 0 & 0 & 0\end{pmatrix},B=\begin{pmatrix}-px_1 & -p\beta x_1 & 0 & -p\beta x_1 & 0\\
0 & 0 & 0 & 0 & 0\\
0 & -p\beta x_3 & -px_3 & 0 & 0\\
0 & 0 & 0 & 0 & 0\\
0 & 0 & 0 & -p\beta x_5 & -px_5\end{pmatrix},$$
$$C=\begin{pmatrix}\frac{cr_1}{r_1+\alpha r_2+\alpha r_4} & 0 & 0 & 0 & 0 \\
\frac{c\alpha r_2}{r_1+\alpha r_2+\alpha r_4} & \frac{cr_2}{r_2+\alpha r_3} & \frac{c\alpha r_2}{\alpha r_2+r_3} & 0 & 0 \\
0 & \frac{c\alpha r_3}{r_2+\alpha r_3} & \frac{cr_3}{\alpha r_2+r_3} & 0 & 0\\
\frac{c\alpha r_4}{r_1+\alpha r_2+\alpha r_4} & 0 & 0 & \frac{cr_4}{r_4+\alpha r_5} & \frac{c\alpha r_4}{\alpha r_4+r_5}\\
0 & 0 & 0 & \frac{c\alpha r_5}{r_4+\alpha r_5} & \frac{cr_5}{\alpha r_4+r_5}\end{pmatrix}=\begin{pmatrix}\frac{br_1}{x_1} & 0 & 0 & 0 & 0\\
\frac{b\alpha r_2}{x_1} & c & c & 0 & 0\\
0 & 0 & 0 & 0 & 0\\
\frac{b\alpha r_4}{x_1} & 0 & 0 & c & c\\
0 & 0 & 0 & 0 & 0\end{pmatrix},$$
$$D=\begin{pmatrix}D_1\\D_2\end{pmatrix},$$
$$D_1=\begin{pmatrix}\frac{cx_1\alpha(r_2+r_4)}{(r_1+\alpha r_2+\alpha r_4)^2}-b & -\frac{cx_1r_1\alpha}{(r_1+\alpha r_2+\alpha r_4)^2} & 0 & -\frac{cx_1r_1\alpha}{(r_1+\alpha r_2+\alpha r_4)^2} & 0\\
-\frac{cx_1\alpha r_2}{(r_1+\alpha r_2+\alpha r_4)^2} & \frac{cx_1\alpha(r_1+\alpha r_4)}{(r_1+\alpha r_2+\alpha r_4)^2}+\frac{cx_2\alpha r_3}{(r_2+\alpha r_3)^2}+\frac{cx_3\alpha r_3}{(\alpha r_2+r_3)^2}-b & -\frac{cx_2r_2\alpha}{(r_2+\alpha r_3)^2}-\frac{cx_3\alpha r_2}{(\alpha r_2+r_3)^2} & -\frac{cx_1\alpha^2r_2}{(r_1+\alpha r_2+\alpha r_4)^2} & 0\\
0 & -\frac{cx_2\alpha r_3}{(r_2+\alpha r_3)^2}-\frac{cx_3r_3\alpha}{(\alpha r_2+r_3)^2} & \frac{cx_2\alpha r_2}{(r_2+\alpha r_3)^2}+\frac{cx_3\alpha r_2}{(\alpha r_2+r_3)^2}-b & 0 & 0\end{pmatrix}$$
$$=\begin{pmatrix}-\frac{b^2r_1}{cx_1} & -\frac{b^2\alpha r_1}{cx_1} & 0 & -\frac{b^2\alpha r_1}{cx_1} & 0\\
-\frac{b^2\alpha r_2}{cx_1} & \frac{b^2\alpha(r_1+\alpha r_4)}{cx_1}-b & b-\frac{b}{\alpha} & -\frac{b^2\alpha^2r_2}{cx_1} & 0\\
0 & 0 & \frac{b}{\alpha}-2b & 0 & 0\\
\end{pmatrix},$$
$$D_2=\begin{pmatrix}
-\frac{c\alpha x_1r_4}{(r_1+\alpha r_2+\alpha r_4)^2} & -\frac{c\alpha^2x_1r_4}{(r_1+\alpha r_2+\alpha r_4)^2} & 0 & \frac{cx_1\alpha(r_1+\alpha r_2)}{(r_1+\alpha r_2+\alpha r_4)^2}+\frac{cx_4\alpha r_5}{(r_4+\alpha r_5)^2}+\frac{cx_5\alpha r_5}{(\alpha r_4+r_5)^2}-b & -\frac{cx_4\alpha r_4}{(r_4+\alpha r_5)^2}-\frac{cx_5\alpha r_4}{(\alpha r_4+r_5)^2}\\
0 & 0 & 0 & -\frac{cx_4\alpha r_5}{(r_4+\alpha r_5)^2}-\frac{cx_5\alpha r_5}{(\alpha r_4+r_5)^2} & \frac{cx_4\alpha r_4}{(r_4+\alpha r_5)^2}+\frac{cx_5\alpha r_4}{(\alpha r_4+r_5)^2}-b
\end{pmatrix}$$
$$=\begin{pmatrix}-\frac{b^2\alpha r_4}{cx_1} & -\frac{b^2\alpha^2r_4}{cx_1} & 0 & \frac{b^2\alpha(r_1+\alpha r_2)}{cx_1}-b & b-\frac{b}{\alpha}\\
0 & 0 & 0 & 0 & \frac{b}{\alpha}-2b
\end{pmatrix}.$$
Let $\lambda_1=f_2-pr_2=f_2-f_3/\beta,\lambda_2=f_4-pr_4=f_4-f_5/\beta,\lambda_3=b/\alpha-2b$, then
$$\det(J-\lambda I)=(\lambda_1-\lambda)(\lambda_2-\lambda)(\lambda_3-\lambda)^2T(\lambda),$$
where 
\begin{align*}T(\lambda)=(\frac{b^2r_1}{cx_1}\lambda+pbr_1)[\lambda^2+b(1-\alpha)\lambda+cp\beta x_5][\lambda^2+b(1-\alpha)\lambda+cp\beta x_3]\\
+\lambda^3[\lambda+b(1-\alpha)][\lambda^2+b(1-\alpha)\lambda+(\frac{b\alpha}{c}\lambda+p\beta x_1)\frac{b\alpha(r_2+r_4)}{x_1}]\\
+cp\beta x_5\lambda^2[\lambda^2+b(1-\alpha)\lambda+\frac{b\alpha r_2}{x_1}(\frac{b\alpha}{c}\lambda+p\beta x_1)]\\
+cp\beta x_3\lambda^2[\lambda^2+b(1-\alpha)\lambda+\frac{b\alpha r_4}{x_1}(\frac{b\alpha}{c}\lambda+p\beta x_1)+cp\beta x_5]
\end{align*}
\begin{align*}
=\lambda^6+[\frac{b^2r_1}{cx_1}(1-\alpha)+b(2-\alpha)]\lambda^5
+b\{f_1+(1-\alpha)[\frac{b^2r_1}{cx_1}(2-\alpha)+b]\}\lambda^4\\
+b^2(1-\alpha)\{2f_1-f_3-f_5+\frac{b}{cx_1}[r_1(b(1-\alpha)+f_3+f_5)+2\alpha^2f_3r_4]\}\lambda^3\\
+b^2(1-\alpha)\{\frac{b^2r_1}{cx_1}(1-\alpha)(f_3+f_5)+pr_1[f_3+f_5+b(1-\alpha)]+f_3f_5(1+\alpha)\}\lambda^2
\end{align*}
\begin{align*}
+b^3(1-\alpha)^2[\frac{br_1}{cx_1}f_3f_5+pr_1(f_3+f_5)]\lambda+pb^3r_1f_3f_5(1-\alpha)^2.
\end{align*}
Detailed computation of $T(\lambda)$ can be found in \ref{tapp}. One can see that all the coefficients are positve, thus $T(\lambda)$ does not have real positive roots. Indeed we can easily find various groups of parameters for which our two LIs stably coexist. For instance, among them are the following two groups. 
\begin{enumerate}
\item $f_1=3,f_2=2,f_3=1,f_4=2,f_5=1,b=1,\alpha=2/3,\beta=4/9$; $\lambda_1=-1/4=\lambda_2<0,\lambda_3=-1/2<0$, $T(\lambda)$ has 3 pairs of complex roots, all with negative real parts.
\item $f_1=4,f_2=1,f_3=2,f_4=1,f_5=1,b=2,\alpha=3/4,\beta=9/16$; $\lambda_1=-23/9<0,\lambda_2=-7/9<0,\lambda_3=-4/3<0$, $T(\lambda)$ has 3 pairs of complex roots, all with negative real parts.
\end{enumerate}
By continuity there are positive measure sets in the parameter space where the LIs coexist stably.

\section{Discussion}\label{dis}

In this paper we proved that local immunodeficiency discovered in \cite{pnas} is a stable and robust phenomenon which may appear already in CRNs with just three types of viruses. Therefore LI should be likely present in all diseases which demonstrate cross-immunodeficiency. It is not necessary to have large CRNs which are typical for Hepatisis C \cite{pnas}.
We also rigorously demonstrated that it is easy to build larger networks with several elements (persistent nodes) which remain invisible to the host's immune system because of their positions in the CRN.

We also demonstrate that LI is a much more general phenomenon than assumed in \cite{pnas}. Indeed a CRN doesn't need to be scale-free \cite{pnas} to produce LI; it just needs a sufficiently non-homogeneous topology. Since our results are built on exact computations for small networks, they leave a little doubt about the presence of stable and robust LI in large CRNs with heterogeneous topology of a general type.

Observe that the phenomenon of local immunodeficiency formally requires the presence of only persistent antigens which manage to escape immune response. However, in all cases with stable and robust LI, altruistic nodes were always present. It is consistent with extensive numerical simulations with large CRNs in \cite{pnas}. Therefore it seems that altruistic antigens are necessary for LI to be a stable and robust phenomenon.

Overall local immunodeficiency seems to be an ubiquitous phenomenon which likely will be present in all diseases demonstrating cross-immunoreactivity.  It calls for future numerical, analytic and, first of all, biological studies. The most important and interesting question is  which types of viruses can play a role of persistent and/or altruistic ones.

\appendix

\section{Computation for symmetric size 2 CRN}\label{2sym}
Consider the symmetric size 2 CRN in Fig.~\ref{2symap}.

\begin {figure}[h]
\begin {tikzpicture}[-latex ,auto ,node distance =2 cm and 2cm ,on grid ,
semithick ,
state/.style ={ circle ,top color =white , bottom color = processblue!20 ,
draw,processblue , text=blue , minimum width =0.5 cm}]
\node[state] (A) {$1$};
\node[state] (B) [right =of A] {$2$};
\path (A) edge [bend left =15] node[above] {} (B);
\path (B) edge [bend left =15] node[] {} (A);
\end{tikzpicture}
\caption{size 2 CRN (symmetric)}
\label{2symap}
\end{figure}

The dynamics of this CRN is described by
$$\begin{cases}
\dot x_1=f_1x_1-px_1(r_1+\beta r_2),\\
\dot x_2=f_2x_2-px_2(\beta r_1+r_2),\\
\dot r_1=c(x_1\frac{r_1}{r_1+\alpha r_2}+x_2\frac{\alpha r_1}{\alpha r_1+r_2})-br_1,\\
\dot r_2=c(x_1\frac{\alpha r_2}{r_1+\alpha r_2}+x_2\frac{r_2}{\alpha r_1+r_2})-br_2.
\end{cases}$$
Consider the fixed point with local immunodeficiency $x_1>0,r_1=0,x_2=0,r_2>0$. One can solve it to be $$x_1=\frac{bf_1}{c\beta},r_1=0,x_2=0,r_2=\frac{f_1}{\beta}.$$
The Jacobian of the system is $$J=\begin{pmatrix}f_1-p(r_1+\beta r_2) & 0 & -px_1 & -p\beta x_1\\
0 & f_2-p(\beta r_1+r_2) & -p\beta x_2 & -px_2\\
\frac{cr_1}{r_1+\alpha r_2} & \frac{c\alpha r_1}{\alpha r_1+r_2} & cx_1\frac{\alpha r_2}{(r_1+\alpha r_2)^2}+cx_2\frac{\alpha r_2}{(\alpha r_1+r_2)^2}-b & -\frac{cx_1\alpha r_1}{(r_1+\alpha r_2)^2}-\frac{cx_2\alpha r_1}{(\alpha r_1+r_2)^2}\\
\frac{c\alpha r_2}{r_1+\alpha r_2} & \frac{cr_2}{\alpha r_1+r_2} & -\frac{cx_1\alpha r_2}{(r_1+\alpha r_2)^2}-\frac{cx_2\alpha r_2}{(\alpha r_1+r_2)^2} & cx_1\frac{\alpha r_1}{(r_1+\alpha r_2)^2}+cx_2\frac{\alpha r_1}{(\alpha r_1+r_2)^2}-b\end{pmatrix}$$
$$=\begin{pmatrix}0 & 0 & -px_1 & -p\beta x_1\\
0 & f_2-pr_2 & 0 & 0\\
0 & 0 & \frac{b}{\alpha}-b & 0\\
c & c & -\frac{b}{\alpha} & -b\end{pmatrix}.$$
$\lambda=\frac{b}{\alpha}-b>0$ is an eigenvalue, so the fixed point is unstable.

\section{Computation for 3-cycle CRN}\label{3nodeap}
The last size three CRN we consider here for illustration is the 3-cycle network in Fig.~\ref{3-cycle}. 

\begin{figure}[h]
\begin{tikzpicture}[-latex ,auto ,node distance =2 cm and 1.15cm ,on grid ,
semithick ,
state/.style ={ circle ,top color =white , bottom color = processblue!20 ,
draw,processblue , text=blue , minimum width =0.5 cm}]
\node[state] (C)
{$3$};
\node[state] (A) [above left=of C] {$1$};
\node[state] (B) [above right =of C] {$2$};
\path (A) edge [] node[left] {} (B);
\path (C) edge [] node[below]{}(A);
\path (B) edge [] node[above] {} (C);
\end{tikzpicture}
\caption{3-cycle CRN}
\label{3-cycle}
\end{figure}

The governing equations in this case are
$$\begin{cases}
\dot x_1=f_1x_1-px_1(r_1+\beta r_2),\\
\dot x_2=f_2x_2-px_2(r_2+\beta r_3),\\
\dot x_3=f_3x_3-px_3(r_3+\beta r_1),\\
\dot r_1=c(x_1\frac{r_1}{r_1+\alpha r_2}+x_3\frac{\alpha r_1}{\alpha r_1+r_3})-br_1,\\
\dot r_2=c(x_1\frac{\alpha r_2}{r_1+\alpha r_2}+x_2\frac{r_2}{r_2+\alpha r_3})-br_2,\\
\dot r_3=c(x_2\frac{\alpha r_3}{r_2+\alpha r_3}+x_3\frac{r_3}{\alpha r_1+r_3})-br_3.
\end{cases}$$

The fixed points of interest are

$$x_1=0,x_2=\frac{bf_2}{cp},x_3=\frac{bf_3}{cp\beta},r_1=\frac{f_3}{p\beta},r_2=\frac{f_2}{p},r_3=0;$$
$$x_1=\frac{bf_1}{cp\beta},x_2=0,x_3=\frac{bf_3}{cp},r_1=0,r_2=\frac{f_1}{p\beta},r_3=\frac{f_3}{p};$$
$$x_1=\frac{bf_1}{cp},x_2=\frac{bf_2}{cp\beta},x_3=0,r_1=\frac{f_1}{p},r_2=0,r_3=\frac{f_2}{p\beta}.$$

The Jacobian of the system equals 
\begin{align*}
J=\begin{pmatrix}f_1-p(r_1+\beta r_2) & 0 & 0 & -px_1 & -p\beta x_1 & 0\\
0 & f_2-p(r_2+\beta r_3) & 0 & 0 & -px_2 & -p\beta x_2\\
0 & 0 & f_3-p(r_3+\beta r_1) & -p\beta x_3 & 0 & -px_3\\
\frac{cr_1}{r_1+\alpha r_2} & 0 & \frac{c\alpha r_1}{\alpha r_1+r_3} & A-b & -\frac{cx_1\alpha r_1}{(r_1+\alpha r_2)^2} & -\frac{cx_3\alpha r_1}{(\alpha r_1+r_3)^2}\\
\frac{c\alpha r_2}{r_1+\alpha r_2} & \frac{cr_2}{r_2+\alpha r_3} & 0 & -\frac{cx_1\alpha r_2}{(r_1+\alpha r_2)^2} & B-b & -\frac{cx_2\alpha r_2}{(r_2+\alpha r_3)^2}\\
0 & \frac{c\alpha r_3}{r_2+\alpha r_3} & \frac{cr_3}{\alpha r_1+r_3} & -\frac{cx_3\alpha r_3}{(\alpha r_1+r_3)^2} & -\frac{cx_2\alpha r_3}{(r_2+\alpha r_3)^2} & C-b\end{pmatrix},\\
 \text{where }A=\frac{cx_1\alpha r_2}{(r_1+\alpha r_2)^2}+\frac{cx_3\alpha r_3}{(\alpha r_1+r_3)^2},B=\frac{cx_1\alpha r_1}{(r_1+\alpha r_2)^2}+\frac{cx_2\alpha r_3}{(r_2+\alpha r_3)^2},C=\frac{cx_2\alpha r_2}{(r_2+\alpha r_3)^2}+\frac{cx_3\alpha r_1}{(\alpha r_1+r_3)^2}.
\end{align*}

At the fixed point $x_1=0,x_2=\frac{bf_2}{cp},x_3=\frac{bf_3}{cp\beta},r_1=\frac{f_3}{p\beta},r_2=\frac{f_2}{p},r_3=0$, we have 
\begin{align*}
A=B=0,C=\alpha b+\frac{b}{\alpha},\\
J=\begin{pmatrix}f_1-\beta f_2-\frac{f_3}{\beta} & 0 & 0 & 0 & 0 & 0\\
0 & 0 & 0 & 0 & -\frac{b}{c}f_2 & -\frac{b\beta}{c}f_2\\
0 & 0 & 0 & -\frac{b}{c}f_3 & 0 & -\frac{b}{c\beta}f_3\\
\frac{cf_3}{f_3+\alpha\beta f_2} & 0 & c & -b & 0 & -\frac{b}{\alpha}\\
\frac{c\alpha\beta f_2}{f_3+\alpha\beta f_2} & c & 0 & 0 & -b & -\alpha b\\
0 & 0 & 0 & 0 & 0 & \frac{b}{\alpha}+\alpha b-b
\end{pmatrix}.
\end{align*}

Because $\lambda=\alpha b+\frac{b}{\alpha}-b>0$ is an eigenvalue this point is unstable.

At the fixed point $x_1=\frac{bf_1}{cp\beta},x_2=0,x_3=\frac{bf_3}{cp},r_1=0,r_2=\frac{f_1}{p\beta},r_3=\frac{f_3}{p}$ we obtain
\begin{align*}
A=\frac{b}{\alpha}+\alpha b,B=C=0,\\
J=\begin{pmatrix}0 & 0 & 0 & -\frac{bf_1}{c\beta} & -\frac{b}{c}f_1 & 0\\
0 & f_2-\frac{f_1}{\beta}-\beta f_3 & 0 & 0 & 0 & 0\\
0 & 0 & 0 & -\frac{b\beta}{c}f_3 & 0 & -\frac{b}{c}f_3\\
0 & 0 & 0 & \frac{b}{\alpha}+\alpha b-b & 0 & 0\\
c & \frac{cf_1}{f_1+\alpha\beta f_3} & 0 & -\frac{b}{\alpha} & -b & 0\\
0 & \frac{c\alpha\beta f_3}{f_1+\alpha\beta f_3} & c & -\alpha b & 0 & -b
\end{pmatrix}.
\end{align*}

Again $\lambda=\frac{b}{\alpha}+\alpha b-b>0$ is an eigenvalue, and this fixed point is unstable.

At the fixed point $x_1=\frac{bf_1}{cp},x_2=\frac{bf_2}{cp\beta},x_3=0,r_1=\frac{f_1}{p},r_2=0,r_3=\frac{f_2}{p\beta}$ we get analogously 
\begin{align*}
A=0,B=\alpha b+\frac{b}{\alpha},C=0,\\
J=\begin{pmatrix}0 & 0 & 0 & -\frac{b}{c}f_1 & -\frac{b\beta}{c}f_1 & 0\\
0 & 0 & 0 & 0 & -\frac{bf_2}{c\beta} & -\frac{b}{c}f_2\\
0 & 0 & f_3-\beta f_1-\frac{f_2}{\beta} & 0 & 0 & 0\\
c & 0 & \frac{c\alpha\beta f_1}{\alpha\beta f_1+f_2} & -b & -\alpha b & 0\\
0 & 0 & 0 & 0 & \alpha b+\frac{b}{\alpha}-b & 0\\
0 & c & \frac{cf_2}{\alpha\beta f_1+f_2} & 0 & -\frac{b}{\alpha} & -b
\end{pmatrix}.
\end{align*}

This fixed point is also unstable because $\lambda=\alpha b+\frac{b}{\alpha}-b>0$ is an eigenvalue.

It is not surprising that for a cyclic network there is no stable local immunodeficiency because this network is invariant with respect to rotations. Therefore it is a  homogeneous network while the networks with local immunodeficiency are characterized by a strong non-homogeneity \cite{pnas}.  
\section{A complete list of fixed points for size 2 and 3 CRNs}\label{listapp}
\begin{itemize}
    \item size 2 CRN
    
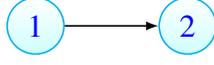
\begin {center}
    \begin {tikzpicture}[-latex ,auto ,node distance =2 cm and 2cm ,on grid,semithick ,
state/.style ={ circle ,top color =white , bottom color = processblue!20 ,
draw,processblue , text=blue , minimum width =0.5 cm}]
\node[state] (A) {$1$};
\node[state] (B) [right =of A] {$2$};
\path (A) edge [] node[above] {} (B);
\end{tikzpicture}
\end{center}
Fixed points:
\begin{enumerate}
\item $$x_1=0,x_2=\frac{bf_2}{cp},r_1=0,r_2=\frac{f_2}{p}$$
\item $$x_1=\frac{bf_1}{cp\beta},x_2=0,r_1=0,r_2=\frac{f_1}{p\beta}$$
\item $$x_1=\frac{bf_1}{cp},x_2=0,r_1=\frac{f_1}{p},r_2=0$$
\item $$f_1=\beta f_2,0<x_1<\frac{bf_2}{cp},x_2=\frac{bf_2}{cp}-x_1,r_1=0,r_2=\frac{f_2}{p}$$
\item $$f_1>\beta f_2,x_1=\frac{b}{cp}(f_1+(\alpha-\beta)f_2),x_2=\frac{bf_2}{cp}(1-\alpha),r_1=\frac{f_1-\beta f_2}{p},r_2=\frac{f_2}{p}$$
\end{enumerate}
\item size 3 CRN
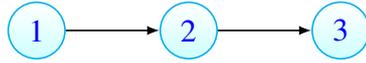
\begin {center}
\begin {tikzpicture}[-latex ,auto ,node distance =02cm and 2cm ,on grid ,
semithick ,
state/.style ={ circle ,top color =white , bottom color = processblue!20 ,
draw,processblue , text=blue , minimum width =0.5 cm}]
\node[state] (A)
{$1$};
\node[state] (B) [right=of A] {$2$};
\node[state] (C) [right =of B] {$3$};
\path (A) edge [] node[below =0.15 cm] {} (B);
\path (B) edge [] node[below =0.15 cm] {} (C);
\end{tikzpicture}
\end{center}
Fixed points:
\begin{enumerate}
\item $$x_1=0,x_2=\frac{bf_2}{cp},x_3=0,r_1=0,r_2=\frac{f_2}{p},r_3=0$$
\item $$f_2>\beta f_3,x_1=0,x_2=\frac{b}{cp}(f_2+(\alpha-\beta)f_3),x_3=\frac{bf_3}{cp}(1-\alpha),r_1=0,r_2=\frac{f_2-\beta f_3}{p},r_3=\frac{f_3}{p}$$
\item $$x_1=\frac{bf_1}{cp\beta},x_2=0,x_3=0,r_1=0,r_2=\frac{f_1}{p\beta},r_3=0$$
\item $$x_1=\frac{bf_1}{cp},x_2=0,x_3=\frac{bf_3}{cp},r_1=\frac{f_1}{p},r_2=0,r_3=\frac{f_3}{p}$$
\item $$x_1=\frac{bf_1}{cp\beta},x_2=0,x_3=\frac{bf_3}{cp},r_1=0,r_2=\frac{f_1}{p\beta},r_3=\frac{f_3}{p}$$
\item $$f_1=\beta f_2,0<x_1<\frac{bf_1}{cp\beta},x_2=\frac{bf_1}{cp\beta}-x_1,x_3=0,r_1=0,r_2=\frac{f_1}{p\beta},r_3=0$$
\item $$x_1=\frac{bf_1}{cp},x_2=\frac{bf_2}{cp\beta},x_3=0,r_1=\frac{f_1}{p},r_2=0,r_3=\frac{f_2}{p\beta}$$
\item $$f_1>\beta f_2,x_1=\frac{b}{cp}(f_1+(\alpha-\beta)f_2),x_2=\frac{bf_2}{cp}(1-\alpha),x_3=0,r_1=\frac{f_1-\beta f_2}{p},r_2=\frac{f_2}{p},r_3=0$$
\item
\begin{align*}
f_1=\beta(f_2-\beta f_3)>0,0<x_1<\frac{b(f_2-\beta f_3)}{cp},x_2=(1+\frac{\alpha f_3}{f_2-\beta f_3})(\frac{b(f_2-\beta f_3)}{cp}-x_1),\\
x_3=\frac{bf_3}{cp}(1-\alpha)+\alpha\frac{f_3}{f_2-\beta f_3}x_1,r_1=0,r_2=\frac{f_2-\beta f_3}{p},r_3=\frac{f_3}{p}
\end{align*}
\item $$f_2=\beta f_3,x_1=\frac{bf_1}{cp},0<x_2<\frac{bf_3}{cp},x_3=\frac{bf_3}{cp}-x_2,r_1=\frac{f_1}{p},r_2=0,r_3=\frac{f_3}{p}$$
\item \begin{align*}f_1>\beta(f_2-\beta f_3)>0,x_1=\frac{b}{cp}(f_1+(\alpha-\beta)(f_2-\beta f_3)),x_2=\frac{b}{cp}(1-\alpha)(f_2+(\alpha-\beta)f_3),\\
x_3=\frac{bf_3}{cp}(1-\alpha(1-\alpha)),r_1=\frac{f_1-\beta f_2+\beta^2 f_3}{p},r_2=\frac{f_2-\beta f_3}{p},r_3=\frac{f_3}{p}
\end{align*}
\end{enumerate}
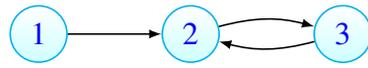
\begin {center}
\begin {tikzpicture}[-latex ,auto ,node distance =2 cm and 2cm ,on grid ,
semithick ,
state/.style ={ circle ,top color =white , bottom color = processblue!20 ,
draw,processblue , text=blue , minimum width =0.5 cm}]
\node[state] (A)
{$1$};
\node[state] (B) [right=of A] {$2$};
\node[state] (C) [right =of B] {$3$};
\path (A) edge [] node[below =0.15 cm] {} (B);
\path (B) edge [bend left =15] node[above] {} (C);
\path (C) edge [bend left =15] node[below] {} (B);
\end{tikzpicture}
\end{center}
Fixed points:
\begin{enumerate}
\item $$x_1=0,x_2=0,x_3=\frac{bf_3}{cp\beta},r_1=0,r_2=\frac{f_3}{p\beta},r_3=0$$
\item $$x_1=0,x_2=\frac{bf_2}{cp},x_3=0,r_1=0,r_2=\frac{f_2}{p},r_3=0$$
\item\begin{align*}
f_3>\beta f_2>\beta^2 f_3,x_1=0,x_2=\frac{b[(1-\alpha\beta)f_2+(\alpha-\beta)f_3]}{cp(1+\alpha)(1-\beta^2)},x_3=\frac{b[(1-\alpha\beta)f_3+(\alpha-\beta)f_2]}{cp(1+\alpha)(1-\beta^2)},\\
r_1=0,r_2=\frac{f_2-\beta f_3}{p(1-\beta^2)},r_3=\frac{f_3-\beta f_2}{p(1-\beta^2)}   
\end{align*}
\item$$f_3=\beta f_2,x_1=0,0<x_2<\frac{bf_2}{cp},x_3=\frac{bf_2}{cp}-x_2,r_1=0,r_2=\frac{f_2}{p},r_3=0$$
\item $$x_1=\frac{bf_1}{cp\beta},x_2=0,x_3=0,r_1=0,r_2=\frac{f_1}{p\beta},r_3=0$$
\item $$f_3=f_1,0<x_1<\frac{bf_1}{cp\beta},x_2=0,x_3=\frac{bf_1}{cp\beta}-x_1,r_1=0,r_2=\frac{f_1}{p\beta},r_3=0$$
\item $$f_3>f_1,x_1=\frac{bf_1}{cp\beta}(1-\alpha),x_2=0,x_3=\frac{b}{cp}(f_3-f_1+\frac{\alpha}{\beta}f_1),r_1=0,r_2=\frac{f_1}{p\beta},r_3=\frac{f_3-f_1}{p}$$
\item $$f_3<f_1,x_1=\frac{b}{cp}(f_1-f_3+\frac{\alpha}{\beta}f_3),x_2=0,x_3=\frac{bf_3}{cp\beta}(1-\alpha),r_1=\frac{f_1-f_3}{p},r_2=\frac{f_3}{p\beta},r_3=0$$
\item $$x_1=\frac{bf_1}{cp},x_2=0,x_3=\frac{bf_3}{cp},r_1=\frac{f_1}{p},r_2=0,r_3=\frac{f_3}{p}$$
\item $$f_1=\beta f_2,0<x_1<\frac{bf_2}{cp},x_2=\frac{bf_2}{cp}-x_1,x_3=0,r_1=0,r_2=\frac{f_2}{p},r_3=0$$
\item $$f_1>\beta f_2,x_1=\frac{b}{cp}(f_1+(\alpha-\beta)f_2),x_2=\frac{bf_2}{cp}(1-\alpha),x_3=0,r_1=\frac{f_1-\beta f_2}{p},r_2=\frac{f_2}{p},r_3=0$$
\item $$x_1=\frac{bf_1}{cp},x_2=\frac{bf_2}{cp\beta},x_3=0,r_1=\frac{f_1}{p},r_2=0,r_3=\frac{f_2}{p\beta}$$
\item $$f_1=f_3=\beta f_2,0<x_1<\frac{bf_2}{cp},0<x_2<\frac{bf_2}{cp}-x_1,x_3=\frac{bf_2}{cp}-x_1-x_2,r_1=0,r_2=\frac{f_2}{p},r_3=0$$
\item
\begin{align*}
(1-\beta^2)f_1=\beta(f_2-\beta f_3)>0,f_3>\beta f_2,0<x_1<b\min\{1-\alpha,\frac{f_2+f_3}{cp(1+\beta)}\},\\x_2=\frac{(1-\alpha\beta)f_2+(\alpha-\beta)f_3}{cp(1+\alpha)(1-\beta^2)}(b-\frac{x_1}{1-\alpha}),
x_3=\frac{(1-\alpha\beta)f_3+(\alpha-\beta)f_2}{cp(1+\alpha)(1-\beta^2)}(b-\frac{\alpha x_1}{1-\alpha}),\\r_1=0,r_2=\frac{f_2-\beta f_3}{p(1-\beta^2)},r_3=\frac{f_3-\beta f_2}{p(1-\beta^2)}
\end{align*}
\item $$f_2=\beta f_3,x_1=\frac{bf_1}{cp},0<x_2<\frac{bf_3}{cp},x_3=\frac{bf_3}{cp}-x_2,r_1=\frac{f_1}{p},r_2=0,r_3=\frac{f_3}{p}$$
\item $$\hspace{-8pt}f_1>\beta f_2=f_3,x_1=\frac{b}{cp}(f_1+(\alpha-\beta)f_2),0<x_2<\frac{bf_2}{cp}(1-\alpha),x_3=\frac{bf_2}{cp}-x_2,r_1=\frac{f_1-\beta f_2}{p},r_2=\frac{f_2}{p},r_3=0$$
\item
\begin{align*}
(1-\beta^2)f_1>\beta(f_2-\beta f_3)>0,f_3>\beta f_2,x_1=\frac{bf_1}{cp}+\frac{b(\alpha-\beta)}{cp(1-\beta^2)}(f_2-\beta f_3),\\
x_2=\frac{b(1-2\alpha)}{cp(1-\alpha^2)(1-\beta^2)}((1-\alpha\beta)f_2+(\alpha-\beta)f_3),\\
x_3=\frac{b(1-\alpha+\alpha^2)}{cp(1-\alpha^2)(1-\beta^2)}((1-\alpha\beta)f_3+(\alpha-\beta)f_2),r_1=\frac{f_1}{p}-\beta\frac{f_2-\beta f_3}{p(1-\beta^2)},\\
r_2=\frac{f_2-\beta f_3}{p(1-\beta^2)},r_3=\frac{f_3-\beta f_2}{p(1-\beta^2)}
\end{align*}
\end{enumerate}

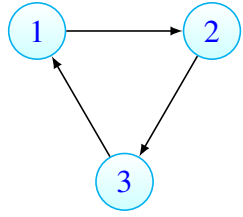
\begin{center}
\begin{tikzpicture}[-latex ,auto ,node distance =2 cm and 1.15cm ,on grid ,
semithick ,
state/.style ={ circle ,top color =white , bottom color = processblue!20 ,
draw,processblue , text=blue , minimum width =0.5 cm}]
\node[state] (C)
{$3$};
\node[state] (A) [above left=of C] {$1$};
\node[state] (B) [above right =of C] {$2$};
\path (A) edge [] node[left] {} (B);
\path (C) edge [] node[below]{}(A);
\path (B) edge [] node[above] {} (C);
\end{tikzpicture}
\end{center}
Fixed points:
\begin{enumerate}
\item $$f_2>\beta f_3,x_1=0,x_2=\frac{b}{cp}(f_2+(\alpha-\beta)f_3),x_3=\frac{bf_3}{cp}(1-\alpha),r_1=0,r_2=\frac{f_2-\beta f_3}{p},r_3=\frac{f_3}{p}$$
\item $$x_1=0,x_2=\frac{bf_2}{cp},x_3=\frac{bf_3}{cp\beta},r_1=\frac{f_3}{p\beta},r_2=\frac{f_2}{p},r_3=0$$
\item $$x_1=\frac{bf_1}{cp\beta},x_2=0,x_3=\frac{bf_3}{cp},r_1=0,r_2=\frac{f_1}{p\beta},r_3=\frac{f_3}{p}$$
\item $$f_3>\beta f_1,x_1=\frac{bf_1}{cp}(1-\alpha),x_2=0,x_3=\frac{b}{cp}(f_3+(\alpha-\beta)f_1),r_1=\frac{f_1}{p},r_2=0,r_3=\frac{f_3-\beta f_1}{p}$$
\item $$x_1=\frac{bf_1}{cp},x_2=\frac{bf_2}{cp\beta},x_3=0,r_1=\frac{f_1}{p},r_2=0,r_3=\frac{f_2}{p\beta}$$
\item $$f_1>\beta f_2,x_1=\frac{b}{cp}(f_1+(\alpha-\beta)f_2),x_2=\frac{bf_2}{cp}(1-\alpha),x_3=0,r_1=\frac{f_1-\beta f_2}{p},r_2=\frac{f_2}{p},r_3=0$$
\item $$\hspace{-40pt}f_2=\frac{f_1}{\beta}+\beta f_3,0<x_1<\frac{bf_1}{cp\beta},x_2=(\frac{b}{c}-\frac{x_1p\beta}{f_1})\frac{f_1+\alpha\beta f_3}{p\beta},x_3=(\frac{b}{c}(1-\alpha)+\frac{\alpha x_1p\beta}{f_1})\frac{f_3}{p},r_1=0,r_2=\frac{f_1}{p\beta},r_3=\frac{f_3}{p}$$
\item $$\hspace{-70pt}f_3=\beta f_1+\frac{f_2}{\beta},(1-\alpha)\frac{bf_1}{cp}<x_1<\frac{bf_1}{cp},x_2=\frac{f_2}{p\alpha\beta}(\frac{x_1p}{f_1}-(1-\alpha)\frac{b}{c}),x_3=(\frac{f_1}{p}+\frac{f_2}{p\alpha\beta})(\frac{b}{c}-\frac{x_1p}{f_1}),r_1=\frac{f_1}{p},r_2=0,r_3=\frac{f_2}{p\beta}$$
\item $$\hspace{-25pt}f_1=\beta f_2+\frac{f_3}{\beta},0<x_1<\frac{b}{c}(r_1+\alpha r_2),x_2=(\frac{b}{c}-\alpha\frac{x_1}{r_1+\alpha r_2})r_2,x_3=r_1(\frac{b}{c}-\frac{x_1}{r_1+\alpha r_2}),r_1=\frac{f_3}{p\beta},r_2=\frac{f_2}{p},r_3=0$$
\item
\begin{align*}
f_1-\beta f_2+\beta^2f_3>0,f_2-\beta f_3+\beta^2 f_1>0,f_3-\beta f_1+\beta^2 f_2>0,\\
x_1=\frac{b}{c(1+\alpha)}(r_1+\alpha r_2),x_2=\frac{b}{c(1+\alpha)}(r_2+\alpha r_3),x_3=\frac{b}{c(1+\alpha)}(r_3+\alpha r_1),\\
r_1=\frac{f_1-\beta f_2+\beta^2f_3}{p(1+\beta^3)},r_2=\frac{f_2-\beta f_3+\beta^2f_1}{p(1+\beta^3)},r_3=\frac{f_3-\beta f_1+\beta^2f_2}{p(1+\beta^3)}
\end{align*}
\end{enumerate}
\end{itemize}

\section{Size 4 mildly asymmetric networks: existence \& stability of LI}
The CRN we consider here is the "T-shaped" network with four nodes in Fig.~\ref{T4}.
\begin {figure}[h]
\begin {tikzpicture}[-latex ,auto ,node distance =2 cm and 2cm ,on grid ,
semithick ,
state/.style ={ circle ,top color =white , bottom color = processblue!20 ,
draw,processblue , text=blue , minimum width =0.5 cm}]
\node[state] (A)
{$1$};
\node[state] (B) [right=of A] {$2$};
\node[state] (C) [below =of A] {$3$};
\node[state] (D) [left =of A] {$4$};
\path (B) edge [] node[] {} (A);
\path (C) edge [] node[] {} (A);
\path (D) edge [] node[] {} (A);
\end{tikzpicture}
\caption{size 4 CRN}
\label{T4}
\end{figure}

For this specific size 4 CRN, we want node 1 to be altruistic, i.e. $x_1=0,r_1>0$. Observe that the nodes 2, 3 and 4 are situated symmetrically. Without loss of generality we may assume that the node 2 is persistent while the nodes 3, 4 are neutral active, i.e. $x_2>0,r_2=0,x_3>0,r_3>0,x_4>0,x_4>0$.

The dynamical equations \eqref{population} assume the form
$$\begin{cases}
\dot x_1=f_1x_1-px_1r_1,\\
\dot x_2=f_2x_2-px_2(\beta r_1+r_2),\\
\dot x_3=f_3x_3-px_3(\beta r_1+r_3),\\
\dot x_4=f_4x_4-px_4(\beta r_1+r_4),\\
\dot r_1=c(x_1+x_2\frac{\alpha r_1}{\alpha r_1+r_2}+x_3\frac{\alpha r_1}{\alpha r_1+r_3}+x_4\frac{\alpha r_1}{\alpha r_1+r_4})-br_1,\\
\dot r_2=cx_2\frac{r_2}{\alpha r_1+r_2}-br_2,\\
\dot r_3=cx_3\frac{r_3}{\alpha r_1+r_3}-br_3,\\
\dot r_4=cx_4\frac{r_4}{\alpha r_1+r_4}-br_4.\\
\end{cases}$$

Under assumptions $f_2<f_3,f_2<f_4,\alpha<1/2$ (so that the population values are positive), we get the fixed point with local immunodeficiency:
$$x_1=0,r_1=\frac{f_2}{p\beta},x_2=\frac{bf_2(1-2\alpha)}{cp\beta},r_2=0,$$
$$x_3=\frac{b}{cp}(\frac{\alpha}{\beta}f_2+f_3-f_2),r_3=\frac{f_3-f_2}{p},x_4=\frac{b}{cp}(\frac{\alpha}{\beta}f_2+f_4-f_2),r_4=\frac{f_4-f_2}{p}.$$
The corresponding Jacobian is,
$$J=\begin{pmatrix}A & B \\ C & D\end{pmatrix},$$
$$A=\begin{pmatrix}f_1-pr_1 & 0 & 0 & 0\\
0 & f_2-p(\beta r_1+r_2) & 0 & 0\\
0 & 0 & f_3-p(\beta r_1+r_3) & 0\\
0 & 0 & 0 & f_4-p(\beta r_1+r_4)\end{pmatrix}=\begin{pmatrix}f_1-pr_1 & 0 & 0 & 0\\
0 & 0 & 0 & 0\\
0 & 0 & 0 & 0\\
0 & 0 & 0 & 0\end{pmatrix},$$
$$B=\begin{pmatrix}-px_1 & 0 & 0 & 0\\
-p\beta x_2 & -px_2 & 0 & 0\\
-p\beta x_3 & 0 & -px_3 & 0\\
-p\beta x_4 & 0 & 0 & -px_4\end{pmatrix}=\begin{pmatrix}0 & 0 & 0 & 0\\
-p\beta x_2 & -px_2 & 0 & 0\\
-p\beta x_3 & 0 & -px_3 & 0\\
-p\beta x_4 & 0 & 0 & -px_4\end{pmatrix},$$
$$C=\begin{pmatrix}c & \frac{c\alpha r_1}{\alpha r_1+r_2} & \frac{c\alpha r_1}{\alpha r_1+r_3} & \frac{c\alpha r_1}{\alpha r_1+r_4}\\
0 & \frac{cr_2}{\alpha r_1+r_2} & 0 & 0\\
0 & 0 & \frac{cr_3}{\alpha r_1+r_3} & 0\\
0 & 0 & 0 & \frac{cr_4}{\alpha r_1+r_4}\end{pmatrix}=\begin{pmatrix}c & c & \frac{b\alpha r_1}{x_3} & \frac{b\alpha r_1}{x_4}\\
0 & 0 & 0 & 0\\
0 & 0 & \frac{br_3}{x_3} & 0\\
0 & 0 & 0 & \frac{br_4}{x_4}\end{pmatrix},$$
$$D=\begin{pmatrix}\frac{c\alpha x_2r_2}{(\alpha r_1+r_2)^2}+\frac{c\alpha x_3r_3}{(\alpha r_1+r_3)^2}+\frac{c\alpha x_4r_4}{(\alpha r_1+r_4)^2}-b & -\frac{c\alpha x_2r_1}{(\alpha r_1+r_2)^2} & -\frac{c\alpha x_3r_1}{(\alpha r_1+r_3)^2} & -\frac{c\alpha x_4r_1}{(\alpha r_1+r_4)^2}\\
-\frac{c\alpha x_2r_2}{(\alpha r_1+r_2)^2} & \frac{c\alpha x_2r_1}{(\alpha r_1+r_2)^2}-b & 0 & 0\\
-\frac{c\alpha x_3r_3}{(\alpha r_1+r_3)^2} & 0 & \frac{c\alpha x_3r_1}{(\alpha r_1+r_3)^2}-b & 0\\
-\frac{c\alpha x_4r_4}{(\alpha r_1+r_4)^2} & 0 & 0 & \frac{c\alpha x_4r_1}{(\alpha r_1+r_4)^2}-b
\end{pmatrix}$$
$$=\begin{pmatrix}\frac{\alpha b^2r_3}{cx_3}+\frac{\alpha b^2r_4}{cx_4}-b & 2b-\frac{b}{\alpha} & -\frac{\alpha b^2r_1}{cx_3} & -\frac{\alpha b^2r_1}{cx_4}\\
0 & \frac{b}{\alpha}-3b & 0 & 0\\
-\frac{\alpha b^2r_3}{cx_3} & 0 & \frac{\alpha b^2r_1}{cx_3}-b & 0\\
-\frac{\alpha b^2r_4}{cx_4} & 0 & 0 & \frac{\alpha b^2r_1}{cx_4}-b
\end{pmatrix}.$$

 As an exact numerical example with a stable local immunodeficiency consider the system's parameters assuming  the following values $b=c=p=1,\alpha=2/5,\beta=4/25,f_1=f_2=1,f_3=f_4=2$. One can compute the Jacobian numerically and see all the eigenvalues are either real negative or complex with negative real parts. It follows by continuity that there exists a positive measure set in the parameter space where this local immunodeficiency is stable.

\section{Detailed computation of $T(\lambda)$}\label{tapp}

After column reduction, we get 
\begin{align*}
T(\lambda)=\begin{vmatrix}-\lambda & 0 & 0 & -px_1-\frac{b}{c}\lambda & -p\beta x_1-\frac{\alpha b}{c}\lambda & -p\beta x_1-\frac{\alpha b}{c}\lambda\\
0 & -\lambda & 0 & 0 & -p\beta x_3 & 0\\
0 & 0 & -\lambda & 0 & 0 & -p\beta x_5\\
\frac{br_1}{x_1} & 0 & 0 & -\lambda & 0 & 0 \\
\frac{b\alpha r_2}{x_1} & c & 0 & 0 & \alpha b-b-\lambda & 0 \\
\frac{b\alpha x_4}{x_1} & 0 & c & 0 & 0 & \alpha b-b-\lambda\end{vmatrix}.
\end{align*}
There are many zeros among these entries. Expanding along the rows or columns with the most number of 0s is the simplest way to compute the determinant. The following computation uses the expansion along the row that has the lowest index number among all rows and columns with the most number of 0s.
\begin{align*}
T(\lambda)=-\lambda\begin{vmatrix}-\lambda & 0 & -px_1-\frac{b}{c}\lambda & -p\beta x_1-\frac{\alpha b}{c}\lambda & -p\beta x_1-\frac{\alpha b}{c}\lambda\\
0 & -\lambda & 0 & 0 & -p\beta x_5\\
\frac{br_1}{x_1} & 0 & -\lambda & 0 & 0 \\
\frac{b\alpha r_2}{x_1} & 0 & 0 & \alpha b-b-\lambda & 0\\
\frac{b\alpha r_4}{x_4} & c & 0 & 0 & \alpha b-b-\lambda\end{vmatrix}\\
+p\beta x_3\begin{vmatrix}
-\lambda & 0 & 0 & -px_1-\frac{b}{c}\lambda & -p\beta x_1-\frac{\alpha b}{c}\lambda\\
0 & 0 & -\lambda & 0 & -p\beta x_5\\
\frac{br_1}{x_1} & 0 & 0 & -\lambda & 0\\
\frac{b\alpha r_2}{x_1} & c & 0 & 0 & 0\\
\frac{b\alpha r_4}{x_1} & 0 & c & 0 & \alpha b-b-\lambda
\end{vmatrix}
\end{align*}

\begin{align*}
=\lambda^2\begin{vmatrix}-\lambda & -px_1-\frac{b}{c}\lambda & -p\beta x_1-\frac{\alpha b}{c}\lambda & -p\beta x_1-\frac{b\alpha}{c}\lambda\\
\frac{br_1}{x_1} & -\lambda & 0 & 0 \\
\frac{b\alpha r_2}{x_1} & 0 & \alpha b-b-\lambda & 0\\
\frac{b\alpha r_4}{x_1} & 0 & 0 & \alpha b-b-\lambda\end{vmatrix}\\
-p\beta x_5\lambda\begin{vmatrix}
-\lambda & 0 & -px_1-\frac{b}{c}\lambda & p\beta _1-\frac{b\alpha}{c}\lambda\\
\frac{br_1}{x_1} & 0 & -\lambda & 0\\
\frac{b\alpha r_2}{x_1} & 0 & 0 & \alpha b-b-\lambda\\
\frac{b\alpha r_4}{x_1} & c & 0 & 0
\end{vmatrix}+cp\beta x_3\begin{vmatrix}
-\lambda & 0 & -px_1-\frac{b}{c}\lambda & -p\beta x_1-\frac{b\alpha}{c}\lambda\\
0 & -\lambda & 0 & -p\beta x_5\\
\frac{br_1}{x_1} & 0 & -\lambda & 0\\
\frac{b\alpha r_4}{x_1} & c & 0 & \alpha b-b-\lambda
\end{vmatrix}
\end{align*}
\begin{align*}
=\lambda^2[-\frac{br_1}{x_1}\begin{vmatrix}-px_1-\frac{b}{c}\lambda & -p\beta x_1-\frac{b\alpha}{c}\lambda & -p\beta x_1-\frac{b\alpha}{c}\lambda\\
0 & \alpha b-b-\lambda & 0\\
0 & 0 & \alpha b-b-\lambda\end{vmatrix}-\lambda\begin{vmatrix}-\lambda & -p\beta x_1-\frac{b\alpha}{c}\lambda & -p\beta x_1-\frac{b\alpha}{c}\lambda\\
\frac{b\alpha r_2}{x_1} & \alpha b-b-\lambda & 0\\
\frac{b\alpha r_4}{x_1} & 0 & \alpha b-b-\lambda\end{vmatrix}]\\
-cp\beta x_5\lambda\begin{vmatrix}-\lambda & -px_1-\frac{b}{c}\lambda & -p\beta x_1-\frac{b\alpha}{c}\lambda\\
\frac{br_1}{x_1} & -\lambda & 0\\
\frac{b\alpha r_2}{x_1} & 0 & \alpha b-b-\lambda\end{vmatrix}
-cp\beta x_3\lambda\begin{vmatrix}
-\lambda & -px_1-\frac{b}{c}\lambda & -p\beta x_1-\frac{b\alpha}{c}\lambda\\
\frac{br_1}{x_1} & -\lambda & 0 \\
\frac{b\alpha r_4}{x_1} & 0 & \alpha b-b-\lambda\end{vmatrix}\\
-cp^2\beta^2x_3x_5\begin{vmatrix}
-\lambda & 0 & -px_1-\frac{b}{c}\lambda\\
\frac{br_1}{x_1} & 0 & -\lambda\\
\frac{b\alpha r_4}{x_1} & c & 0
\end{vmatrix}
\end{align*}
\begin{align*}
=\lambda^2\{\frac{br_1}{x_1}(px_1+\frac{b}{c}\lambda)(\lambda+b-\alpha b)^2\\
+\lambda[\frac{b\alpha r_2}{x_1}(p\beta x_1+\frac{b\alpha}{c}\lambda)(\lambda+b-\alpha b)+(\lambda+b-\alpha b)(\lambda^2+(b-\alpha b)\lambda+(p\beta x_1+\frac{b\alpha}{c}\lambda)\frac{b\alpha r_4}{x_1})]\}\\
+cp\beta x_5\lambda\{\frac{br_1}{x_1}(px_1+\frac{b}{c}\lambda)(\lambda+b-\alpha b)+\lambda[\lambda^2+(b-\alpha b)\lambda+(p\beta x_1+\frac{b\alpha}{c}\lambda)\frac{b\alpha r_2}{x_1}]\}\\
+cp\beta x_3\lambda\{\frac{br_1}{x_1}(px_1+\frac{b}{c}\lambda)(\lambda+b-\alpha b)+\lambda[\lambda^2+(b-\alpha b)\lambda+(p\beta x_1+\frac{b\alpha}{c}\lambda)\frac{b\alpha r_4}{x_1}]\}\\
+c^2p^2\beta^2x_3x_5(\lambda^2+\frac{b^2r_1}{cx_1}\lambda+pbr_1)
\end{align*}
\begin{align*}
=\frac{br_1}{x_1}(px_1+\frac{b}{c}\lambda)(\lambda+b-\alpha b)[\lambda^2(\lambda+b-\alpha b)+cp\beta x_5\lambda+cp\beta x_3\lambda]\\
+(p\beta x_1+\frac{b\alpha}{c}\lambda)[\frac{b\alpha r_2}{x_1}\lambda^3(\lambda+b-\alpha b)+\frac{b\alpha r_4}{x_1}\lambda^3(\lambda+b-\alpha b)+\frac{b\alpha r_2}{x_1}cp\beta x_5\lambda^2+\frac{b\alpha r_4}{x_1}cp\beta x_3\lambda^2]\\
+\lambda^4(\lambda+b-\alpha b)^2+cp\beta (x_5+x_3)\lambda^3(\lambda+b-\alpha b)+c^2p^2\beta^2x_3x_5[\lambda^2+\frac{b^2r_1}{cx_1}\lambda+pbr_1]
\end{align*}
\begin{align*}
=\frac{br_1}{x_1}\lambda(\frac{b}{c}\lambda+px_1)(\lambda+b-\alpha b)[\lambda^2+(b-\alpha b)\lambda+cp\beta(x_3+x_5)]\\
+\lambda^2(\frac{b\alpha}{c}\lambda+p\beta x_1)[(\frac{b\alpha r_2}{x_1}+\frac{b\alpha r_4}{x_1})\lambda(\lambda+b-\alpha b)+\frac{b\alpha}{x_1}cp\beta(r_2x_5+r_4x_3)]\\
+\lambda^3(\lambda+b-\alpha b)[\lambda^2+(b-\alpha b)\lambda+cp\beta(x_3+x_5)]\\
+c^2p^2\beta^2x_3x_5(\lambda^2+\frac{b^2r_1}{cx_1}\lambda+pbr_1)
\end{align*}
\begin{align*}
=\lambda[\lambda+b(1-\alpha)][\lambda^2+b(1-\alpha)\lambda+b(1-\alpha)(f_3+f_5)][\lambda^2+\frac{b^2r_1}{cx_1}\lambda+pbr_1]\\
+\lambda^2(\frac{b\alpha}{c}\lambda+p\beta x_1)[(c-\frac{br_1}{x_1})\lambda(\lambda+b-\alpha b)+2\frac{b^2\alpha(1-\alpha)}{x_1}r_4f_3]\\
+b^2(1-\alpha)^2f_3f_5(\lambda^2+\frac{b^2r_1}{cx_1}\lambda+pbr_1).\\
\end{align*}

\section*{Acknowledgements}
We are indebted to P. Skums and Yu. Khudyakov for valuable discussions and suggestions. We also want to thank both anonymous referees for valuable questions and comments. This work was partially supported by the NSF grant CCF-BSF-1664836 and by the NIH grant 1R01EB025022.




\begin{thebibliography}{9}

\bibitem{pnas}
P. Skums, L. Bunimovich, Y. Khudyakov, Antigenic cooperation among intrahost HCV variants organized into a complex network of cross-immunoreactivity, Proceedings of the National Academy of Sciences of the United States of America 112 (21) (2015) 6653-6658.


\bibitem{hattori}
M. Hattori, K. Yashioka, T. Aiyama, K. Iwata, Y. Terazawa, M. Ishigami, M. Yano, S. Kakumu, Broadly reactive antibodies to hypervariable region 1 in hepatitis C virus-infected patient sera: relation to viral loads and response to interferon, Hepatology 27 (6) (1998) 1703-1710.

\bibitem{campo}
D. S. Campo, Z. Dimitrova, L. Yamasaki, P. Skums, D. T. Lau, G. Vaughan, J. C. Forbi, C.-G. Teo, Y. Khudyakov, Next-generation sequencing reveals large connected networks of intra-host HCV variants, BMC Genomics 15 (Suppl 5) (2014) S4.


\bibitem{nowak00}
M. A. Nowak, R. M. May, Virus Dynamics: Mathematical Principles of Immunology and Virology, Oxford University Press, 2000.

\bibitem{yoshioka}
K. Yoshioka, T. Aiyama, A. Okumura, M. Takayanagi, K. Iwata, T. Ishikawa, Y. Nagai, S. Kakumu, Humoral Immune Response to the Hypervariable Region of Hepatitis C Virus Differs between Genotypes 1b and 2a, The Journal of Infectious Diseases 175 (3) (1997) 505-510.

\bibitem{wodarz}
D. Wodarz, Hepatitis C virus dynamics and pathology: the role of CTL and antibody responses, Journal of General Virology 84 (Pt 7) (2013) 1743-1750.

\bibitem{nowak90}
M. A. Nowak, R. M. May, R. M. Anderson, The evolutionary dynamics of HIV-1 quasispecies and the development of immunodeficiency disease, AIDS 4 (11) (1990) 1095-1103.

\bibitem{nowak91}
M. A. Nowak, R. M. May, Mathematical biology of HIV infections: antigenic variation and diversity threshold, Mathematical Biosciences, 106 (1) (1991) 1-21.

\bibitem{nowak91s}
M. A. Nowak, R. M. Anderson, A. R. Mclean, R. M. May, Antigenic Diversity Thresholds and the Development of AIDS, Science 254 (5034) (1991) 963-969.

\bibitem{tom}
T. Francis, Jr., On the Doctrine of Original Antigenic Sin, Proceedings of the American Philosophical Society 104 (6) (1960) 572-578.

\bibitem{pan}
K. Pan, Understanding Original Antigenic Sin in Influenza with a Dynamical System, PLoS ONE 6 (8) (2011) e23910.

\bibitem{shin}
B. Rehermann, E.-C. Shin, Private aspects of heterologous immunity, Journal of Experimental Medicine 201 (5) (2005) 667-670.

\bibitem{hvi}
M. S. Parsons, S. Muller, H. Kholer, M. D. Grant, N. F. Bernard, On the benefits of sin:Can greater understanding of the 1F7-idiotypic repertoire freeze enhance HIV vaccine development?, Human Vaccines and Immunotherapeutics 9 (7) (2013) 1532-1538.

\bibitem{jin}
J. H. Kim, I. Skountzou, R. Compans, J. Jacob, Original Antigenic Sin Responses to Influenza Viruses, J Immunol 183 (5) (2009) 3294-3301.

\bibitem{alotofauthors}
C. M. Midgley, M. Bajwa-Joseph, S. Vasanawathana, W. Limpitikul, B. Wills, A. Flanagan, E. Waiyaiya, H. B. Tran, A. E. Cowper, P. Chotiyarnwon, J. M. Grimes, S. Yoksan, P. Malasit, C. P. Simmons, J. Mongkolsapaya, G. R. Screaton, An In-Depth Analysis of Original Antigenic Sin in Dengue Virus Infection, Journal of Virology 85 (1) (2011) 410-421.


\end{thebibliography}

\end{document}